\documentclass[12pt]{article}
\usepackage{hyperref,amsfonts,amssymb,amsmath}
\usepackage{mathtools}

\usepackage{tikz-cd}
\usepackage{tikz-cd, amsmath, amssymb}

\usepackage{color}

\def\tr{\mathrm{tr}}

\usepackage{soul}

\def\C{\mathbb{C}}
\def\F{\cal{F}}

\def\Z{\mathbb{Z}}
\def\Q{\mathbb{Q}}

\def\R{\mathbb{R}}
\def\P{\mathbb{P}}
\def\A{\mathbb{A}}

\def\I{{{\mathcal I}}}

\def\bq{ \begin{equation} }
\def\eq{ \end{equation} }
\def\ben{ \begin{eqnarray} }
\def\en{ \end{eqnarray} }

\def\frac#1#2{{#1\over #2}}

\def\on#1#2{\mathop{\vbox{\ialign{##\crcr\noalign{\kern2pt}
$\scriptstyle{#2}$\crcr\noalign{\kern2pt\nointerlineskip}
\kern-2pt$\hfil\displaystyle{#1}\hfil$\crcr}}}\limits}

\begin{document}

\title{Multiplication kernels}
\author{Maxim Kontsevich,  Alexander Odesskii}
   \date{}
\vspace{-20mm}
   \maketitle
\vspace{-7mm}
\begin{center}
IHES, 35 route de Chartres, Bures-sur-Yvette, F-91440,
France  \\[1ex]
and \\[1ex]
Brock University, 1812 Sir Isaac Brock Way, St. Catharines, ON, L2S 3A1 Canada\\[1ex]
e-mails: \\
\texttt{maxim@ihes.fr}\\
\texttt{aodesski@brocku.ca}
\end{center}

\medskip

\begin{abstract}

We introduce the notion of multiplication kernels of birational and $D$-module type and give various examples. We also introduce the notion of a semi-classical multiplication kernel associated with an integrable system and discuss its quantization. Finally, we discuss geometric and algebraic aspects of method of separation of 
variables, and describe hypothetically a cyclic $D$-module for the generalized multiplication kernels for Hitchin systems for groups $GL_r$.

\end{abstract}

\newpage

\tableofcontents

\newpage

 \section{Introduction}
 
  \subsection{Informal explanation of the problem}
  
  Suppose we have a collection of commuting linear operators\footnote{A lot of interesting examples of such collections appears in the theory of integrable systems and in representation theory.} $T_{\alpha}$ acting on a finite dimensional vector space $V$. Assume that the joint spectrum of these operators is simple. Then we have a basis of $V$ (of joint eigenvectors of $T_{\alpha}$) defined up to permutation and rescaling. Assume furthermore that we choose a vector $v\in V$ which is cyclic, i.e. generates $V$ as the module over the algebra generated by $T_{\alpha}$. Then the basis of eigenvectors $e_i$ is defined only up to permutations if we put the constraint $v=\sum_i e_i$ i.e. $v$ has coordinates $(1,1,...,1)$. Alternatively, one may assume that a cyclic covector $u\in V^*$ is given 
  and normalize basic elements by the condition $(u,e_i)=1$. 
  
  The basis $\{e_i\}$ up to permutations can be encoded by the structure of a {\it commutative associative unital algebra} on $V$ with the multiplication $e_i\cdot e_j=\delta_{i,j}e_i$ where $\delta_{i,j}$ is the Kronecker delta. 
  
  In this paper we deal with a functional analogue of this situation where vector space $V$ is a space of functions (in a broad sense) in one or several variables, or 
  more generally, $V$ is a space of functions on a smooth or algebraic manifold. The typical example is $V=C^{\infty}(\R^n)$ where the commuting operators are 
  derivations $\frac{\partial}{\partial x_k},~k=1,...,n$. The continuous analogue of joint eigenvectors consists of Fourier
  modes $e_{\lambda}(x)=e^{ix\cdot \lambda},~\lambda\in \R^n$ where the normalization is given by $e_{\lambda}(0)=1$. The multiplication is given by 
  $e_{\lambda}(x) * e_{\mu}(x)=\delta(\lambda-\mu)e_{\lambda}(x)$ where $\delta(\lambda-\mu)$ is the Dirac delta\footnote{We denote our multiplication in functional case by $*$ in order to distinguish it from usual pointwise product of functions.}. Notice that this multiplication is the additive convolution and can be written in terms of the standard multiplication of functions as 
  $$f * g(y)=\int_{\R^n\times \R^n} \delta(x_1+x_2-y)f(x_1)g(x_2)dx_1dx_2.$$
  
  Let $\F$ be a vector space of functions on a manifold $X$. We denote a typical element of $\F$ by $f(x)\in \F$, where $x\in X$. 
  
  We want to study 
  commutative associative multiplications $\F\otimes \F\to\F$ on the vector space $\F$. We write such multiplication in the form 
  \begin{equation} \label{mult}
  f*g(y)=\int_{X\times X} K(x_1,x_2,y) f(x_1)g(x_2) dx_1dx_2
  \end{equation}
  where $K(x_1,x_2,y)$ is a kernel of our multiplication and $dx$ is a measure on $X$ defined by a volume form. Commutativity of our multiplication means
  \begin{equation} \label{com}
  K(x_1,x_2,y)=K(x_2,x_1,y)
  \end{equation}
 and associativity of our multiplication means
  \begin{equation} \label{ass}
 \int_X K(x_1,x_2,y)K(y,x_3,z)dy=\int_X K(x_1,x_3,y)K(y,x_2,z)dy
  \end{equation} 
  In which way one can verify the associativity condition (\ref{ass}) if the kernel $K$ is given in a closed form (for example, in terms of elementary functions)? In principle, one can suggest the following possibilities:
  
  {\bf 1.} Compute explicitly the l.h.s. and the r.h.s. of the equation (\ref{ass}). 
  
  {\bf 2.} Find a change of variables $y\mapsto \tilde{y}$ which transforms the l.h.s. of the equation 
  $$K(x_1,x_2,y)K(y,x_3,z)~dy=K(x_1,x_3,\tilde{y})K(\tilde{y},x_2,z)~d\tilde{y}$$
  to its r.h.s. In this case the equation (\ref{ass}) also holds provided that $X$ is a cycle.
  
  {\bf 3.} Prove that the l.h.s and the r.h.s. of the equation (\ref{ass}) satisfies the same holonomic system of differential equations as a function in $x_1,x_2,x_3,z$. Strictly speaking, this does not mean that the equation (\ref{ass}) holds on the nose, but still can be considered as an associativity condition for a kernel $K$.
  
  In this paper we do not deal with analysis as the possibility 1 suggests, and concentrate on algebraic side 
  of the problem. For example, assuming that $X$ is a small circle in $\C$ around zero, when exploring possibility 1 we can  reduce the integration to an algebraic operations if $K(x_1,x_2,y)\in\C[\frac{1}{y}][[x_1,x_2]]$, i.e. $K$ is a power series in $x_1,x_2$ with coefficients polynomial in $\frac{1}{y}$.  
  
  Exploring the possibility 2 we assume that $X$ is an algebraic variety and a change of variables $y\mapsto \tilde{y}$ defines a birational mapping $X\to X$. 
  
  Finally, in the possibility 3 we assume that $X$ is a cycle and, therefore $\int_X \frac{\partial h}{\partial q} dq=0$ for any $h$. This assumption reduces computations to algebraic manipulations in differential algebra.
  
  More generally, the kernel $K$ can be given by integration of an ``elementary function'' over some auxiliary variables.   Namely, let $Q$ be another manifold with a measure $dq$ defined by a volume form. Assume that
 \begin{equation} \label{kerint}
 K(x_1,x_2,y)=\int_Q K(x_1,x_2,y,q)dq.
 \end{equation}
 In this case associativity condition takes a form
\begin{equation} \label{ass2}
 \underset{X\times Q\times Q}{\int} K(x_1,x_2,y,q_1)K(y,x_3,z,q_2)dydq_1dq_2=\underset{X\times Q\times Q}{\int} K(x_1,x_3,y,q_3)K(y,x_2,z,q_4)dydq_3dq_4.
  \end{equation}

 Finally, let $X$, $Q$ both be algebraic varieties over $\C$. In this case we assume that 
 \begin{equation} \label{elfun}
 K=K_1^{s_1}...K_l^{s_l}
  \end{equation} 
 where $s_1,...,s_l\in\C$ are arbitrary parameters and $K_1,...,K_l$ are either algebraic functions or exponential of algebraic functions\footnote{In the case of general ground field we replace $K_i^{s_i}$ by $\chi_i(K_i)$ where $\chi_i$ are either multiplicative or additive characters of the ground field.}. These lead to the following
  
  {\bf Definition 1.1.1.} We say that $K(x_1,x_2,y,q)$ is a multiplication kernel of birational type if there exists a birational automorphism $(y,q_1,q_2)\to (\tilde{y},\tilde{q}_3,\tilde{q}_4)$ of $X\times Q\times Q$ which transforms the l.h.s. of the equation
  $$K(x_1,x_2,y,q_1)K(y,x_3,z,q_2)dydq_1dq_2= K(x_1,x_3,\tilde{y},\tilde{q}_3)K(\tilde{y},x_2,z,\tilde{q}_4)d\tilde{y}d\tilde{q}_3d\tilde{q}_4$$
  to its r.h.s. If $K=K_1^{s_1}...K_l^{s_l}$, then this birational automorphism should not depend on $s_1,...,s_l$.
  
  The kernel $K(x_1,x_2,y)$ given by (\ref{kerint}),   (\ref{elfun}) satisfies a holonomic system of differential equations in $x_1,x_2,y$ i.e. it gives a holonomic $D$-module endowed with a cyclic vector. The associativity constraint (\ref{ass}) can be understood as an isomorphism between two holonomic $D$-modules in $x_1,x_2,x_3,y$ with cyclic vectors. The operation of integration corresponds to the direct image of $D$-modules.   
  
  {\bf Definition 1.1.2.} In the situation as above, we say that $K(x_1,x_2,y)$ is a multiplication kernel of $D$-module type.
  
 {\bf Remark 1.1.1.} The associativity conditions from Definition 1.1.2 can be formulated less abstractly as follows. Assume that $K(x_1,x_2,y)$ is a solution of a holonomic system of differential equations\footnote{We understand integral in (\ref{kerint}) as direct image which means that our system of differential equations for $K(x_1,x_2,y)$ is a consequence of a holonomic  system of differential equations for $K(x_1,x_2,y,q)$.} in $x_1,x_2,y$. Let $I^{12}$ be the left ideal in the ring of differential operators in $x_1,x_2,y,x_3,z$ generated by differential equations of $K(x_1,x_2,y)$. Let $I^{23}$ be the left ideal in the same ring generated by differential equations of $K(y,x_3,z)$. It is clear that $I^{23}$ is obtained from $I^{12}$ by the change of variables $(x_1,x_2,y)\to (y,x_3,z)$. Define a left ideal $I^{123}$ in the ring of differential operators in $x_1,x_2,x_3,z$ by
  $$I^{123}=(I^{12}+I^{23}+J)\cap R_{123}$$
  where $J$ is the {\it right} ideal generated by $\frac{\partial}{\partial y}$ and $R_{123}$ is the ring of differential operators in $x_1,x_2,x_3,z$. By construction, $I^{123}$ consists of differential equations for 
  $$\int_{X\times Q\times Q} K(x_1,x_2,y,q_1)K(y,x_3,z,q_2)dydq_1dq_2$$
  where $X\times Q\times Q$ is a cycle.
  
 {\bf Definition 1.1.2$^{\prime}$.}  We say that $K(x_1,x_2,y)$ is a multiplication kernel of $D$-module type if $I^{123}$ is invariant with respect to interchanging 
  of $x_2$ and $x_3$.

   \subsection{Formal setups for functions and integrations}
   
   In the previous informal definitions of multiplication kernel we use the notions of function, integration etc in non-rigorous way. 
   There are various rigorous formalisms for the notions of ``explicit formulas'' and ``function'' for an  algebraic variety $X$ defined over a field $k$, none of which is totally satisfactory. 
   
   {\bf 1)} If $k$ is a local field (i.e. $k=\R,\C$, or a finite extension of $\Q_p$ or $\mathbb{F}_p(\!(t)\!)$) and $\pi:~Y\to X$ is a family of $n$-dimensional varieties over $X$, 
   endowed with the volume element $vol\in\Gamma(Y, K_{Y/X})$ along fibers of $\pi$. Here $K_{Y/X}=\Lambda^n T^*_{Y/X}$ is the relative canonical line bundle. Then we obtain an $\R$-valued function on the set $X(k)$ given by 
   $\pi_{*}(|vol|)$, the integration of the density $|vol|$ along fibers. Here we assume that the integral is convergent, at least for the generic point $x\in X(k)$. 
   
   One can twist the integral by additive and multiplicative characters applied to rational functions on $Y$ (this is a formal replacement of exponentials and fractional powers). Moreover, it is enough to assume more generally that for some integer $N\geq 1$ we have an element of $\Gamma(Y,K^{\otimes N}_{Y/X})$, e.g. $vol^{\otimes N}$ 
   if $vol$ is defined up to multiplication by $N$-th root of 1. 
   
   This definition has several drawbacks. First, the same function on $X(k)$ can be presented as an integral in different ways, and it is not clear whether the equality always follows from a birational equivalence. Second, one expects that some interesting functions (for example in representation theory) do not have such an integral 
   representation (see \cite{kazd} for discussion). Third, there is a problem with $S_n$-covariance for $n$-fold compositions of multiplication kernels, see 
   Section 2.6 for details. 
   
   {\bf 2)} If $char(k)=0$, we can encode a ``function'' by a holonomic $D$-module endowed with a cyclic vector, see Section 2.7 for details. The drawback of this definition is that we can not distinguish functions $f$ and $c f$ where $c$ is a non-zero constant. 
   
   {\bf 3)} We can forget the cyclic vector in the previous definition, and encode a function by an {\it equivalence class} of  holonomic $D$-modules. 
   
   ${\bf 3^{\prime}}${\bf)} The previous definition can be transported to arbitrary characteristic, with holonomic $D$-modules replaced by motivic constructible sheaves. 
   
   {\bf 4)} Finally, if we are in the case of positive characteristic $p>0$ and $k=\mathbb{F}_{p^n}$, then  one can associate with motivic constructible sheaf a 
   $\overline{\Q}^{CM}$ valued function\footnote{Here $\overline{\Q}^{CM}\subset\overline{\Q}$ stands for the maximal totally real extension of $\Q$ with added 
   $\sqrt{-1}$. This notation comes from theory of complex multiplication for abelian varieties.} on the finite set $X(\mathbb{F}_{p^n})$ given by the trace of Frobenius.
   Surprisingly, here we have again a well-defined function (as in {\bf 1)}) although the information on the cyclic vector seems to be lost. 
   
   In all formalisms above one can speak about integrals (as direct images) and hence the associativity constraint (\ref{ass}) makes sense. Therefore, one can speak about multiplication kernels in different contexts.

   \subsection{Multiplication formulas for special functions}
   
   Let us discuss a dual viewpoint on multiplication kernels. The product $*$ defined by (\ref{ass}) on the space $\F$ of functions gives by duality a coproduct $\Delta$ on 
   the dual space  $\F^*$ of densities. The continuous basis $e_{\lambda}(x)$ of elementary projectors for $*$ gives the dual basis  $e^*_{\lambda}(x) dx$ of $\F^*$. 
   The property $\Delta(e^*_{\lambda}(y) dy)=e^*_{\lambda}(x_1) dx_1\otimes e^*_{\lambda}(x_2) dx_2$ is equivalent to\footnote{In the examples below we make an identification between $\F$ and $\F^*$ and write this formula in terms of $e_{\lambda}$.} 
   $$e^*_{\lambda}(x_1)  e^*_{\lambda}(x_2)=\int K(x_1,x_2,y) e^*_{\lambda}(y) dy$$
   for all $\lambda$, where $K(x_1,x_2,y)$ is the same kernel as in (\ref{ass}) and, in particular, does not depend on $\lambda$.
   
   Let us give several examples where $M$ is one-dimensional and eigenfunctions  $e_{\lambda}(x)$ are classical special functions. In these examples $T_{\alpha}$ 
   consists of one differential operator $T$.
   
   {\bf Example 1.3.1.} $M=\R$,~~~ $T=\frac{d}{dx}$,~~~ $e_{\lambda}(x)=e^{\lambda x}$. We have $Te_{\lambda}(x)=\lambda e_{\lambda}(x)$, the normalization is defined by $e_{\lambda}(0)=1$, and 
   $$e_{\lambda}(x_1)e_{\lambda}(x_2)=e_{\lambda}(x_1+x_2)=\int_{\R} \delta(x_1+x_2-y)e_{\lambda}(y)dy.$$
   
   {\bf Example 1.3.2.} $M=\R_{>0}$,~~~ $T=x\frac{d}{dx}$,~~~ $e_{\lambda}(x)=x^{\lambda}$. We have $Te_{\lambda}(x)=\lambda e_{\lambda}(x)$, the normalization is defined by $e_{\lambda}(1)=1$, and 
   $$e_{\lambda}(x_1)e_{\lambda}(x_2)=e_{\lambda}(x_1x_2)=\int_{\R} \delta(x_1x_2-y)e_{\lambda}(y) dy.$$
   
  {\bf Example 1.3.3.} $M=\R_{>0}$,~~~  $T=\frac{d^2}{dx^2}+x^{-1}\frac{d}{dx}$,~~~ $e_{\lambda}(x)=J_0(\lambda x)$. We have $Te_{\lambda}(x)=-\lambda^2 e_{\lambda}(x)$, the normalization is defined by $\lim_{x\to 0} e_{\lambda}(x)=1$, and 
   $$e_{\lambda}(x_1)e_{\lambda}(x_2)=\int^{x_2+x_1}_{|x_2-x_1|} \frac{e_{\lambda}(y)}{A(x_1,x_2,y)}\frac{ydy}{2\pi}.$$  
   Here 
   $$J_0(x)=\sum_{m=0}^{\infty}\frac{(-1)^m}{m!^2}\Big(\frac{x}{2}\Big)^{2m}=\frac{1}{2\pi i}\oint e^{x\frac{u-u^{-1}}{2}}\frac{du}{u}$$
   is Bessel function and 
   $$A(x_1,x_2,y)=\frac{1}{4}(2x_1^2x_2^2+2x_1^2y^2+2x_2^2y^2-x_1^4-x_2^4-y^4)^{\frac{1}{2}}$$
   is the area of triangle with sides of the length $x_1,x_2,y$. This multiplication formula is called Sonine–Gegenbauer formula \cite{R,SG}.

  \subsection{Relations to Langlands correspondence and integrable systems}
  
  Theory of automorphic forms provides examples of commuting operators (Hecke operators). In the case of a curve $C$ over finite field $\mathbb{F}_q$ these operators act
  on the space of functions on the countable set of isomorphism classes of $G$-bundles on $C$ where $G$ is a reductive group. In the case $G=GL_2$ there are multiplicity 
  one theorems which guarantee that the joint spectrum is simple. In the old paper \cite{K} of one of us, the multiplication kernel was written explicitly in a special 
  case of rank 2 bundles on $\P^1$ with the parabolic structure in 4 points. Also, in the same paper, a kernel for the case of local field $k$ was given by the formula
  $$f*g(y)=\int_{y\in k,F_t(x,y,z)\in (k^*)^2} \frac{f(x_1)g(x_2)}{|F_t(x_1,x_2,y)|^{\frac{1}{2}}}|dy|$$
  where 
  $$F_t(x_1,x_2,y)=(x_1x_2+x_1y+x_2y-t)^2+4x_1x_2y(1+t-x_1-x_2-y).$$
  Recently, in the paper \cite{efk}, Hecke operators in the case of curves over $\C$ were defined (but not yet the multiplication kernels). 
  
  The joint spectrum of commuting integral Hecke operators for the case of curves over local fields is rather mysterious, and its relation to the usual Langlands program 
  is quite unclear. In the case of a non-archimedean field $k$ with the residue field $\mathbb{F}_p$, a finite ``low frequency'' part of the spectrum is presumably the 
  same as the spectrum for the case of curves over finite fields and hence is related to Galois representations. For $k=\C$ the joint spectrum is expected to coincide \cite{efk}
  with the set of opers (roughly speaking, differential equations of rank $r$) with real monodromy. Similarly, for $k=\R$ the joint spectrum is expected to be the 
  spectrum of the algebra of commuting differential operators on the set of $\R$-points of algebraic variety $Bun_G$, coming from the quantization of Hitchin
  integrable system.
  
  In general, if $H_1,...,H_n$ are commuting differential operators on $n$-dimensional manifold, then for any scalar parameters $\lambda_1,...,\lambda_n$ we have a 
  holonomic system 
  $$(H_i-\lambda_i)\psi(x)=0,~i=1,...,n.$$
  Let us denote $\psi_{\lambda}(x)$ a solution of this system where $\lambda=(\lambda_1,...,\lambda_n)$. In order to have unique solution one has to impose some 
  normalization conditions. See Section 4 for details.
  
  In the case $G=GL_r$ and arbitrary constraints at singularities, there exists a remarkable birational symplectomorphism between the phase space of Hitchin 
  integrable systems and the cotangent bundle to $Sym^gC$ where $g$ is the genus of the generic spectral curve or, equivalently, the dimension of the base of the
   integrable 
  system. This construction is called the method of separation of variables \cite{Sk}. It is expected that in the case $G=GL_r$ there exists an integral operator given by 
  an explicit kernel, identifying functions on $Bun_G$ and on $Sym^gC$. Moreover, eigenfunctions of commuting differential operators on $Bun_G$ (or of Hecke operators) 
  map to symmetric functions on $Sym^gC$ of the form $\phi_{\lambda}(x_1)...\phi_{\lambda}(x_g)$ which are external powers of functions in one variable. 
  
  In this presentation the multiplication kernel can be expressed by subsequent integration in terms of a more elementary kernel which we denote by $K_{g+1,g}$, which is a function of $2g+1$ 
  variables whereas the original multiplication kernel is a function of $3g$ variables. We will discuss in details this approach in Section 4.
  
  \subsection{Multiplication kernels in other contexts}
  
  The first group of questions we want to discuss here is related to commuting families of Hecke operators in theory of modular forms. The multiplicity one theorems 
  in the theory of automorphic forms for group $GL_r$ are valid not only for the case of curves over finite field, but also in the number field case. This leads e.g. 
  to the following question concerning classical modular forms for group $SL(2,\Z)$. For any $n\geq 1$ there are $d_n$ Hecke eigenforms 
  $$f_i^{(n)}(q)=q+\sum_{j>1}a_{i,j}^{(n)}q^j$$ 
  of weight $n$ and level 1, where $a_{i,j}^{(n)},~i=1,...,d_n$ are eigenvalues of Hecke operator 
  $H_j^{(n)}$ and $d_n$ is the dimension of the space of cusp forms of weight $n$. Coefficients of these forms are algebraic integers, not necessarily rational.  Let us consider the generating series in 4 variables 
  $$K=\sum_{n\geq 1}t^n~\sum_{i=1}^{d_n}f_i^{(n)}(q_1)f_i^{(n)}(q_2)f_i^{(n)}(q_3)\in \overline{\Q}[[q_1,q_2,q_3,t]].$$
  One can show that $K\in \Z[[q_1,q_2,q_3,t]].$ We can also write $K$ in terms of traces of products of Heche operators as 
  $$K=\sum_{n,j_1,j_2,j_3\geq 1} \tr(H_{j_1}^{(n)}H_{j_2}^{(n)}H_{j_3}^{(n)})~q_1^{j_1}q_2^{j_2}q_3^{j_3}t^n.$$
   It will be interesting to find a closed formula for $K$. A similar question can be asked about higher level modular forms and 
  about Maass forms. 
  
  In order to explain the second group of questions, we start with an example. Let $X=[0,1]$ be the unit interval in $\R$, and $K:X^3\to \{0,1\}$ be the characteristic 
  function of the closed tetrahedron with vertices $(0,0,0),(0,1,1),(1,0,1),(1,1,0)$.
  $$K(x,y,z)=\begin{cases}~ 1,~\text{if}~ x\leq y+z,~y\leq x+z,~ z\leq x+y,~ x+y+z\leq 1 \\~ 0~\text{~otherwise.} \end{cases}$$
  Then the following is true: for any integer $n\geq 1$ introduce vector space $$A_n=\Q^{X\cap \frac{1}{n}\Z}=\Q^{\{0,\frac{1}{n},...,1\}}=\Q^{n+1}$$ 
  with the basis 
  $e_x,x\in X\cap \frac{1}{n}\Z$. Define a product in $A_n$ by 
  $$e_{x_1}\cdot e_{x_2}=\sum_{x_3\in X\cap \frac{1}{n}\Z} K(x_1,x_2,x_3)e_{x_3}. $$
  Then this product is commutative and associative. In fact, it is Verlinde algebra for $sl_2$ at level $n$. 
  
  The proof of associativity can be made independent of $n$. Namely, it follows from the existence of a piecewise-linear identification with integer coefficients 
  ($\Z$PL in short) of two 5-dimensional polytopes fibered over $X^4$. Let
  $$P_1=\{(x_1,x_2,x_3,x_4,y);~(x_1,x_2,y)\in K,~(y,x_3,x_4)\in K\},$$
  $$P_2=\{(x_1,x_2,x_3,x_4,y);~(x_1,x_3,y)\in K,~(y,x_2,x_4)\in K\}.$$
  Define two maps $\pi_i:P_i\to X^4$, $i=1,2$ by $\pi_i(x_1,x_2,x_3,x_4,y)=(x_1,x_2,x_3,x_4)$. One can check that for all $x_1,x_2,x_3,x_4$ the fibers 
  $\pi^{-1}_1(x_1,x_2,x_3,x_4)$ and $\pi^{-1}_2(x_1,x_2,x_3,x_4)$  are closed intervals of the same length, and can be identified by a shift. The resulting map 
  $P_1\to P_2$ is a $\Z$PL homeomorphism.  This argument is similar to the cut-and-paste proof of associativity in the case of multiplication kernels for varieties 
  over finite field studied in \cite{K}. 
  
  This example leads to several questions: 
  
  {\bf 1.} ~Generalize it to the case of other reductive groups, 
  
  {\bf 2.} ~Find other $\Z$PL examples of multiplication kernels, 
  
  {\bf 3.} Find the relation 
  with multiplication kernels given by integral operators over non-archimedian fields, 
  
  {\bf 4.} Find similar formulas for multiplication kernels where numbers of integer points 
  in polytopes is replaced by volumes of polytopes. 
  
   \subsection{Other questions for functional analogue of tensor algebra}
   
   In this paper we study explicit associative commutative kernels using purely algebraic framework (see discussion after associativity condition (\ref{ass})). This  algebraic approach can be applied to other problems in functional analogue of tensor algebra. 
   
   {\bf 1.} Let ${\F}_1$, ${\F}_2$ be two spaces of functions, possibly on different manifolds. One can study explicit kernels for mappings 
   $$R:~{\F}_1\otimes {\F}_2\to {\F}_2$$
   subject to constraint $R(f_1,R(f_2,g))=R(f_2,R(f_1,g))$
   where $f_1,f_2\in {\F}_1,~g\in {\F}_2$. This condition means that all linear operators on ${\F}_2$ of the form $g\mapsto R(f,g)$ commute. Note that associative 
   commutative kernels provide examples of this structure in the case ${\F}_1={\F}_2$ because operators of multiplication by a given element commute in commutative associative algebras.
   
   {\bf 2.} Let $e_{\lambda}(x)$ be the set of joint eigenfunctions of a family of commuting operators on a space of functions $\F$. Here we do not need to choose a normalization.  Define a mapping 
   $$R:~\F\otimes \F\to \F\otimes\F$$
   by $R(e_{\lambda}(x)\otimes e_{\mu}(x))=\delta(\lambda-\mu)e_{\lambda}(x)\otimes e_{\lambda}(x).$ One can study explicit kernels for this mapping for some 
   interesting classes of commuting operators. Notice that $R$ also satisfies to an analogue of associativity and commutativity constraints: $R$ is 
   invariant with respect to the action of $S_2\times S_2$ and the composition $R^{13}\circ R^{23}$ is invariant with respect to the action of $S_3\times S_3$. A hypothetical example of such structure is given in Remark 3.3.2.
   
   {\bf 3.} Given $g\geq 1$ one can study kernels for mapping 
   $$R:~{\F}^{\otimes (g+1)}\to {\F}^{\otimes g}$$
   invariant with respect to the action of $S_{g+1}\times S_g$ such that $R\circ (R\otimes Id_{\F})$ is invariant with respect to the action of $S_{g+2}\times S_g$. 
   We call this structure a {\it generalized product}. Such a product induces a structure of an associative commutative algebra on the space $Sym^g\F$ (see Proposition 
   4.1.1).
   
   {\bf 4.} One can study kernels for associative but not necessarily commutative multiplications. 
   
   {\bf 5.} Let ${\F}_1$, ${\F}_2$ be two spaces of functions. One can study kernels for two inverse linear mappings ${\F}_1\to {\F}_2$ and ${\F}_2\to {\F}_1$. 
   This question can be though of as a fundamental question of integral geometry in the sense of Gelfand-Gindikin-Graev \cite{GGG}.

  \subsection{Content of the paper}
  
  In Section 2 we define the notion of multiplication kernel of birational type in more rigorous way and give examples. The most part of examples are related to the 
  Hitchin systems for group $GL_2$ on the curve $\P^1$ with 4 or more regular singular points. Another way of constructing examples is to solve certain functional 
  equations which is explained in Remark 2.2.2. Notice that Sections 2.2 - 2.5 can be read independently of other parts of the paper provided that the reader is fine 
  with informal explanation of Definition 1.1.1 from Section 1.1.
  
  In Section 3 we explain that classical integrable systems give commutative monoids in the Weinstein category of symplectic varieties and Lagrangian correspondences. We 
  also discuss quantization which is a construction of a multiplication kernel of $D$-module type starting from a quantum intergable system. Our examples are again 
  related with Hitchin systems. It seems to be an interesting and important problem to find explicitly multiplication kernels for various quantum integrable systems. 
  
  In Section 4 we describe semi-classical geometry of Sklyanin's method of separation of variables. The algebraic counterpart of this method is the notion of a 
  generalized product, a map $Sym^{g+1}V\to Sym^g V$ satisfying an analog of associativity constraint. We also introduce a hypotetical construction of quantum generalized 
  products via formal solutions of differential equations expanded in a chosen base point. All considerations in Section 4 can be generalized to trigonometric and 
  elliptic difference equations which is beyond the standard geometric Langlands perspective.
  
  In Section 5 we explain, and illustrate by examples, how to construct multiplication kernels satisfying the property (\ref{ass}) with $M$ is a small circle 
  $|z|=\varepsilon$, $0<\varepsilon \ll 1$ starting with a differential operator. These kernels can also be considered as lifts of more abstract kernels from Section 4. 
  Notice that  Section 5 can be read  independently of other parts of the paper.

\section{Multiplication kernels of birational type}

\subsection{Formulation of the problem}

Let us fix a ground field $k$ of characteristic zero, and a commutative algebraic group\footnote{In our examples $A$ is the product of the additive group scheme $\mathbb{G}_a$ and of several 
copies of the multiplicative group scheme $\mathbb{G}_m$.}  $A$ over $k$. 

Consider the following category ${\cal{C}=C}_{k,A}$. Its objects are smooth equividimensional varieties over $k$. The set of morphisms $Hom_{\cal{C}}(X_1,X_2)$ is 
defined as the set of equivalence classes of tuples $(Z,~\pi_1:Z\to X_1,~\pi_2:Z\to X_2,~ vol,~ \rho:Z\to A)$ where $Z$ is an equividimensional smooth variety over $k$, 
$\pi_1,\pi_2$ are smooth morphisms (submersions), $vol\in\Gamma (Z,K_{Z/X_2})/(\pm 1)$ is a volume element along fibers of $\pi_2$ up to a sign\footnote{More generally, 
one can modify the definition by replacing $vol\in\Gamma(X,K_{Z/X_2})$ by its pover $vol^{\otimes N}\in\Gamma(X,K^{\otimes N}_{Z/X_2})$ for $N\geq 1$.}, 
and $\rho:Z\to A$ is an 
arbitrary map of varieties. The equivalence relation is generated by identifications of tuples 
\begin{equation}   \label{eqrel}
(Z, \pi_1, \pi_2, vol, \rho)\sim (U, \pi_1|_{U}, \pi_2|_U, vol|_U, \rho|_U) 
\end{equation}
for fixed
$X_1,~X_2$ where $U\subset Z$ is a Zariski open dense  subvariety of $Z$. 

{\bf Remark 2.1.1.} If the ground field $k$ is a local field, then the morphisms in $\cal C$ can be thought of as formal integral operators depending on a generic character 
of locally compact abelian group $A(k)$. Namely, for a variety $X/k$ denote by $\F_X$ the space of $\C$-valued continuous functions on $X(k)$. A tuple 
$(Z,~\pi_1:Z\to X_1,~\pi_2:Z\to X_2,~ vol,~ \rho:Z\to A)$ as above gives a formal integral operator $K:\F_{x_1}\to \F_{X_2}$ given by 
$$K(f)(x_2)=\int_{z\in\pi^{-1}_2(x_2)(k)}f(\pi_2(z))\cdot \chi(\rho(z))\cdot |vol|_{\pi_2^{-1}(x_2)}$$
where $\chi:~A(k)\to \C^*$ is a character. We ignore the convergence issues here. This heuristics explains the following definition of the composition in $\cal C$.

For two tuples $$(Z,\pi_1:Z\to X_1,\pi_2:Z\to X_2, vol, \rho:Z\to A),$$ $$(Z^{\prime},\pi^{\prime}_2:Z^{\prime}\to X_2,\pi^{\prime}_3:Z^{\prime}\to X_3, vol^{\prime}, \rho^{\prime}:Z^{\prime}\to A)$$ their composition is given by the fibered product 
$$Z^{\prime\prime}=Z\times_{X_2}Z^{\prime}$$
endowed with maps $\pi_1^{\prime\prime}=\pi_1\circ\pi_2^*(\pi^{\prime}_1):Z^{\prime\prime}\to X_1$, 
$\pi_3^{\prime\prime}=\pi^{\prime}_3\circ(\pi^{\prime}_1)^*(\pi_2):Z^{\prime\prime}\to X_3$. 

This can be represented by the following commuting diagram:
$$
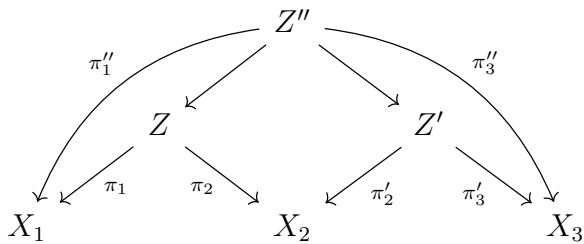
\begin{tikzcd}
  & &  Z''\arrow{rd}\arrow{ld}  \arrow[bend right,lldd,"\pi''_1"']   \arrow[bend left,rrdd,"\pi''_3"]  & & \\
  & Z \arrow[rd,"\pi_2"']\arrow[ld,"\pi_1"] & & Z' \arrow[rd,"\pi_3'"']\arrow[ld,"\pi'_2"] &\\
  X_1 & & X_2 & & X_3
 \end{tikzcd}$$

The volume element $vol^{\prime\prime}$ along fibers of $\pi^{\prime\prime}_3$ 
is obtained by the multiplication of volume elements along fibers of maps $Z^{\prime\prime}\to Z^{\prime}$ and $Z^{\prime}\to X_3$. The map $\chi^{\prime\prime}:Z^{\prime\prime}\to A$ is defined as the product in group scheme $A$ of maps $\rho\circ(Z^{\prime\prime}\to Z)$ and 
$\rho^{\prime}\circ(Z^{\prime\prime}\to Z^{\prime})$.

One can check that the composition is well-defined on equivalence classes of representatives of morphisms. 

The identity morphism $id_{\C}(X)$ is given by $Z=X$, $pr_1=pr_2=id:X\to X$, $vol=1$, and $\rho(x)=0\in A$. 

One can see that two varieties $X_1$, $X_2$ are isomorphic as objects of $\cal C$ iff they are birationally equivalent.

{\bf Remark 2.1.2.} In the definition of ${\cal C}$ one can omit the condition that $\pi_1$ is submersion and the equivalence relation generated by (\ref{eqrel}). In the 
modified category one looses the birational invariance. On the other hand, if we want to keep the equivalence relation (and birational invariance), then the composition 
is defined if we assume that $\pi_1$ is dominant on each component of $Z$. Passing to a Zariski open dense set $U\subset Z$, we can replace this condition by smoothness 
of $\pi_1$. 

Category $\cal C$ carries the natural structure of a symmetric monoidal category . The tensor product on objects is given by the usual product of varieties. In the 
definition of the tensor product of morphisms we use the product in group scheme $A$.

{\bf Definition 2.1.1.} A multiplication kernel of birational type is a commutative semigroup object\footnote{In the framework of category $({\cal C},\otimes)$ it is not reasonable to require an object to be a commutative {\it monoid}, which is a commutative semigroup object with a unit.} in $({\cal C},\otimes)$.

{\bf Remark 2.1.3.}  This definition seems to be too general, as we get some pathological examples. The issue is related to the fact that in our heuristics with local field, we ignore the question of convergence. As a first approximation to a better definition (which takes the convergence into account) one can suggest the following. 

{\bf Definition 2.1.2.} A morphism $(Z,~\pi_1:Z\to X_1,~\pi_2:Z\to X_2,~ vol,~ \rho:Z\to A)$ in $\cal C$ is called {\it geometrically convergent} if the following property holds.
Consider generic point $x_2\in X_2$. Then $\pi_2^{-1}(x_2)$ is a smooth variety which is mapped to $X_1\times A$ by $(\pi_1,\rho)$. Let $Y_{x_2}$ denote the image 
of a connected component of $\pi_2^{-1}(x_2)$. This is not necessarily a smooth variety, but it is nevertheless smooth at its generic point  $y\in Y_{x_2}$. On the
smooth variety $V_{x_2,y}=(\pi_1,\rho)^{-1}(y)\subset Y_{x_2}$ we have a volume element $vol_{x_2,y}$  defined up to 
multiplication by a non-zero constant. Namely, $vol_{x_2,y}$ is defined as the ratio of $vol_{\pi^{-1}_2(x_2)}$ and a non-zero element in $\Lambda^dT^*_yY_{x_2}$ where 
$d=\dim Y_{x_2}$. We demand that $vol_{x_2,y}$ extends to volume form without poles on some (or equivalently, on all) smooth compactification of $V_{x_2,y}$. 

This definition of convergence is not completely satisfactory. For example, the composition of geometrically convergent morphisms is not necessarily geometrically 
convergent. This reflects the fact that in the case of local field and $A=0$ the integral operator associated with a geometrically convergent morphism maps the space 
of bounded measurable functions to a larger space of unbounded measurable functions. One can not compose such operators in general. 

On the other hand, there exist situations when the integral operator corresponding to a not geometrically convergent morphism gives a well-defined compact operator. This 
is related to the fact that integrals $\int_{V(k)}|vol|$ for a meromorphic volume element $vol$ on algebraic variety $V/K$ can be convergent even if $vol$ has poles on 
$V(\bar{k})$. For example, the sphere $S^{n-1}=\{(x_1,...,x_n)\in\R^n;~\sum_{i=1}^n x_i^2=1\}$ has finite volume. 

For a multiplication kernel of birational type $\mu\in Hom_{\cal C}(X\times X,X)$ one can ask that $\mu_2=\mu$, and $\mu_3=\mu\circ (\mu\otimes id_X)\in Hom_{\cal C}(X^3,X),...,\mu_n\in Hom_{\cal C}(X^n,X),...$ are geometrically convergent where $\mu_n$ corresponds to a product of length $n$.

\subsection{One dimensional examples without auxiliary integration}

Let $X=\mathbb{A}^1$, $Z=X^3$ or a ramified covering of $X^3$. In this case we do not have auxiliary integration in the formulas for kernels.

{\bf Example 2.2.1.} Let 
$$K(x_1,x_2,y)=e^{c(x_1x_2y+y)}~\frac{1}{y}$$
where $c$ is a constant. Then we have
$$K(x_1,x_2,y)~K(y,x_3,z)~dy=K(x_1,x_3,\tilde{y})~K(\tilde{y},x_2,z)~d\tilde{y}$$
if $\tilde{y}=\frac{x_1x_2+x_3z+1}{x_1x_3+x_2z+1}y.$

The example 2.2.1 looks degenerate because $\mu_n$ are not geometrically convergent for $n\geq 4$. It looks plausible, however, that further examples in this Section, as 
well as in Sections 2.3, 2.4 are not degenerate in the sense that all $\mu_n,~n\geq 2$ are geometrically convergent. 

{\bf Example 2.2.2.} Here we assume $A=\mathbb{G}_a$ in notations of Section 2.1.  Let 
$$K(x_1,x_2,y)=e^{c\big(x_1x_2y+\frac{x_1}{x_2y}+\frac{x_2}{x_1y}+\frac{y}{x_1x_2}\big)}~\frac{1}{y}$$
where $c$ is a constant. Then we have
$$K(x_1,x_2,y)~K(y,x_3,z)~dy=K(x_1,x_3,\tilde{y})~K(\tilde{y},x_2,z)~d\tilde{y}$$
if $\tilde{y}=\frac{x_1x_2+x_3z}{x_1x_3+x_2z}y.$

A closely related version of this kernel is 
$$K(x_1,x_2,y)=e^{c\frac{x_1x_2y+x_1+x_2+y}{\sqrt{x_1x_2y}}}\frac{1}{y}.$$

{\bf Example 2.2.3.} Here $A=\mathbb{G}_m^2$. Introduce the notation (see also the formula for $F_t$ in Section 1.4)
$$f_t(x,y,z)= (xy+yz+zx-t)^2+4xyz(1+t-(x+y+z))$$
where $t\ne 0,1$ is a parameter.
Let 

$K(x_1,x_2,y)=$
$$\Big(\frac{x_1x_2(2y-1)-(x_1+x_2)y+t+w_{12}}{(x_1-1)(x_2-1)y}\Big)^{c_1}\Big(\frac{x_1x_2(2y-t)-t(x_1+x_2)y+t^2+tw_{12}}{(x_1-t)(x_2-t)y}\Big)^{c_2}\frac{1}{w_{12}}$$
where $w_{12}:=f_t(x_1,x_2,y)^{1/2}$ and $c_1,c_2$ are arbitrary constants. Then 
\begin{equation} \label{ass5}
K(x_1,x_2,y)K(y,x_3,x_4)~dy=K(x_1,x_3,\tilde{y})K(\tilde{y},x_2,x_4)~d\tilde{y}
\end{equation}
where $\tilde{y}$ is a function in $x_1,x_2,x_3,x_4,y$ independent of $c_1,c_2$. 

More precisely, let $E_1\subset \A^3$ be an affine elliptic curve given by\footnote{Here $w_{12},w_{34},y$ are affine coordinates on $\A^3$ and $t,x_1,x_2,x_3,x_4$ are parameters.} 
$$w_{12}^2=f_t(x_1,x_2,y),~~~w_{34}^2=f_t(y,x_3,x_4)$$
and $E_2\subset \A^3$ be an affine elliptic curve given by\footnote{Here $w_{13},w_{24},\tilde{y}$ are affine coordinates on $\A^3$ and $t,x_1,x_2,x_3,x_4$ are parameters.} 
$$w_{13}^2=f_t(x_1,x_3,\tilde{y}),~~~w_{24}^2=f_t(\tilde{y},x_2,x_4).$$
Then there exists a unique birational mapping $\rho:~E_1\to E_2$ which transforms the l.h.s. of (\ref{ass5}) to its r.h.s. In particular, $j$-invariants of the elliptic curves $E_1$, $E_2$ are equal.  This birational mapping has the form
$$\rho:~(y,w_{12},w_{34})\mapsto (\tilde{y},w_{13},w_{24})$$
where $(y,w_{12},w_{34})\in E_1$, $(\tilde{y},w_{13},w_{24})\in E_2$ and $\tilde{y},w_{13},w_{24}$ are given by
$$\tilde{y}=\frac{(x_1-x_2)(x_3-x_4)y^2-w_{12}w_{34}+(x_1x_2-t)(x_3x_4-t)+\frac{y}{(x_1-x_4)(x_2-x_3)}Q_1}{2(x_1-x_3)(x_2-x_4)y+\frac{2(x_1x_2-x_3x_4)}{(x_1-x_4)(x_2-x_3)}Q}$$
$$w_{13}=\frac{(x_3-x_4)w_{12}y-(x_1-x_2)w_{34}y+\frac{w_{12}}{(x_1-x_3)(x_2-x_3)}Q_2+\frac{w_{34}}{(x_1-x_3)(x_1-x_4)}Q_3}{2(x_2-x_4)y+\frac{2(x_1x_2-x_3x_4)}{(x_1-x_3)(x_1-x_4)(x_2-x_3)}Q}$$
$$w_{24}=\frac{(x_1-x_2)w_{34}y-(x_3-x_4)w_{12}y+\frac{w_{34}}{(x_2-x_3)(x_2-x_4)}Q_4+\frac{w_{12}}{(x_2-x_4)(x_1-x_4)}Q_5}{2(x_1-x_3)y+\frac{2(x_1x_2-x_3x_4)}{(x_2-x_3)(x_2-x_4)(x_1-x_4)}Q}.$$
Here $Q,Q_1,Q_2,Q_3,Q_4,Q_5$ are irreducible polynomials in $x_1,x_2,x_3,x_4,t$ defined by the following properties:
$$(x_2-x_4)w_{13}-(x_1-x_3)w_{24}=(x_3-x_4)w_{12}-(x_1-x_2)w_{34},$$
$$\rho:~(0,x_1x_2-t,x_3x_4-t)\mapsto (0,x_1x_3-t,x_2x_4-t),$$
$$\rho:~(1,x_1x_2-x_1-x_2+t,x_3x_4-x_3-x_4+t)\mapsto (1,x_1x_3-x_1-x_3+t,x_2x_4-x_2-x_4+t),$$
$$\rho:~(t,x_1x_2-tx_1-tx_2+t,x_3x_4-tx_3-tx_4+t)\mapsto (t,x_1x_3-tx_1-tx_3+t,x_2x_4-tx_2-tx_4+t)$$
where 
$$(0,x_1x_2-t,x_3x_4-t),~(1,x_1x_2-x_1-x_2+t,x_3x_4-x_3-x_4+t),$$
$$(t,x_1x_2-tx_1-tx_2+t,x_3x_4-tx_3-tx_4+t)\in E_1,$$
$$(0,x_1x_3-t,x_2x_4-t),~(1,x_1x_3-x_1-x_3+t,x_2x_4-x_2-x_4+t),$$
$$(t,x_1x_3-tx_1-tx_3+t,x_2x_4-tx_2-tx_4+t)\in E_2$$
are rational points of the elliptic curves.

{\bf Remark 2.2.1.} The kernel in Example 2.2.3 depends on four points of $\P^1$ which are set to $0,1,t,\infty$. After an arbitrary fractional linear transformation of the variables $x_1,x_2,x_3,x_4,y,\tilde{y}$ we obtain four pairwise distinct arbitrary points in $\P^1$. Colliding some of these points we obtain various degenerations of the kernels in this family, including the kernel in Example 2.2.2. 

{\bf Remark 2.2.2.} Let $$K(x_1,x_2,y)=\phi(x_1,x_2,y)^c~\psi(x_1,x_2,y)$$ where $\phi,\psi$ are symmetric with respect to $x_1,x_2$. Assume that there exists a function \\ $\tilde{y}(x_1,x_2,x_3,y,z)$ independent of $c$ such that 
$$K(x_1,x_2,y)K(y,x_3,z)~dy=K(x_1,x_3,\tilde{y})K(\tilde{y},x_2,z)~d\tilde{y}.$$
This condition gives a system of functional equations for the functions $\phi,\psi,\tilde{y}$. For computational purposes it is convenient to set 
$$\tilde{y}=y+q_1(x_1,x_2,y,z)\cdot (x_2-x_3)+q_2(x_1,x_2,y,z)\cdot(x_2-x_3)^2+...$$
and assume that the equations 
$$\phi(x_1,x_2,y)\phi(y,x_3,z)=\phi(x_1,x_3,\tilde{y})\phi(\tilde{y},x_2,z),$$
$$\psi(x_1,x_2,y)\psi(y,x_3,z)~dy=\psi(x_1,x_3,\tilde{y})\psi(\tilde{y},x_2,z)~d\tilde{y}$$
hold simultaneously.

Examples 2.2.1, 2.2.3 were obtained by solving this system of functional equations, and Examples 2.5.1, 2.5.2 below were obtained by solving the similar functional 
equations for kernels of the form $$K(x_1,x_2,x_3,y)=\phi(x_1,x_2,x_3,y)^c\psi(x_1,x_2,x_3,y).$$ It looks that any solution of these functional equations is either listed
 in Examples 2.2.1, 2.2.3, 2.5.1, 2.5.2 or can be obtained as a limit of the family described in Example 2.2.3. It would be interesting to study similar functional 
 equations for more general kernels.

\subsection{One dimensional example with auxiliary integration}

Let $X=\mathbb{A}^1$, $Z=\mathbb{A}^5$, $A=\mathbb{G}_m^4$. In this case we have auxiliary integration over 2-dimensional domain in the formulas for kernels. 
Fix $t\ne 0,1$ as in Example 4. Let $s_1,s_2,s_3,r$ be arbitrary constants symbolising a generic character of $A$.

{\bf Example 2.3.1.} Let 
$$K(x_1,x_2,y,q_1,q_2)=(x_1x_2)^{1-s_1}((x_1-1)(x_2-1))^{1-s_2}((x_1-t)(x_2-t))^{1-s_3}F(u,v)$$
where 
$$u=\frac{(x_1-1)(x_2-1)(y-1)}{(t-1)^2},~v=\frac{(x_1-t)(x_2-t)(y-t)}{t(t-1)^2}$$
and
$$F(u,v)=(1-q_1-q_2)^{s_1-r-1}q_1^{s_2-r}q_2^{s_3-r}(q_1q_2+vq_1+uq_2)^{r-2}.$$

{\bf Theorem 2.3.1.} There exists a birational mapping 
$$\mu:~(y,q_1,q_2,q_3,q_4)\to (\tilde{y},\tilde{q_1},\tilde{q_2},\tilde{q_3},\tilde{q_4})$$
depending on $x_1,x_2,x_3,z$, which transforms the l.h.s. of the equation 
\begin{equation} \label{ass6}
\begin{split}
K(x_1,x_2,y,q_1,q_2)K(y,x_3,z,q_3,q_4)~dy\wedge dq_1\wedge dq_2\wedge dq_3\wedge dq_4=\\
\pm K(x_1,x_3,\tilde{y},\tilde{q_1},\tilde{q_2})K(\tilde{y},x_2,z,\tilde{q_3},\tilde{q_4})~d\tilde{y}\wedge d\tilde{q_1}\wedge d\tilde{q_2}\wedge d\tilde{q_3}\wedge d\tilde{q_4}
\end{split}
\end{equation}
to its r.h.s.

{\bf Proof.} The l.h.s. of the equation (\ref{ass6}) can be written as 
$$f_1^{s_1}f_2^{s_2}f_3^{s_3}f_4^r\cdot h~ dy\wedge dq_1\wedge dq_2\wedge dq_3\wedge dq_4$$
where
$$f_1=\frac{1-q_1-q_2}{x_1x_2}~\cdot~ \frac{1-q_3-q_4}{yx_3}$$
$$f_2=\frac{q_1}{(x_1-1)(x_2-1)}~\cdot~ \frac{q_3}{(y-1)(x_3-1)}$$
$$f_3=\frac{q_2}{(x_1-t)(x_2-t)}~\cdot~ \frac{q_4}{(y-t)(x_3-t)}$$
$$f_4=\frac{\Big(q_1q_2+\frac{(x_1-t)(x_2-t)(y-t)}{t(t-1)^2}q_1+\frac{(x_1-1)(x_2-1)(y-1)}{(t-1)^2}q_2\Big)\Big(q_3q_4+\frac{(y-t)(x_3-t)(z-t)}{t(t-1)^2}q_3+\frac{(y-1)(x_3-1)(z-1)}{(t-1)^2}q_4\Big)}{(1-q_1-q_2)q_1q_2\cdot (1-q_3-q_4)q_3q_4}$$
$$h=\frac{x_1x_2(x_1-1)(x_2-1)(x_1-t)(x_2-t)}{\Big(1-q_1-q_2\Big)\Big(q_1q_2+\frac{(x_1-t)(x_2-t)(y-t)}{t(t-1)^2}q_1+\frac{(x_1-1)(x_2-1)(y-1)}{(t-1)^2}q_2\Big)^2}~\cdot$$ 
$$\frac{yx_3(y-1)(x_3-1)(y-t)(x_3-t)} {\Big(1-q_3-q_4\Big)\Big(q_3q_4+\frac{(y-t)(x_3-t)(z-t)}{t(t-1)^2}q_3+\frac{(y-1)(x_3-1)(z-1)}{(t-1)^2}q_4\Big)^2}$$
 The r.h.s. of the equation (\ref{ass6}) can be written as 
$$\tilde{f_1}^{s_1}\tilde{f_2}^{s_2}\tilde{f_3}^{s_3}\tilde{f_4}^r\cdot \tilde{h}~ d\tilde{y}\wedge d\tilde{q_1}\wedge d\tilde{q_2}\wedge d\tilde{q_3}\wedge d\tilde{q_4}$$
where $\tilde{f_1},\tilde{f_2},\tilde{f_3},\tilde{f_4}, \tilde{h}$ are obtained from $f_1,f_2,f_3,f_4,h$ by swapping $x_2$ and $x_3$ and replacing $y,q_1,q_2,q_3,q_4$ by $\tilde{y},\tilde{q_1},\tilde{q_2},\tilde{q_3},\tilde{q_4}$. 

Let $E$ be a curve in affine space $\A^5$ with coordinates $y,q_1,q_2,q_3,q_4$ given by  
$$f_1=C_1,~f_2=C_2,~f_3=C_3,~f_4=C_4.$$

Let $\tilde{E}$ be a curve in affine space $\A^5$ with coordinates $\tilde{y},\tilde{q_1},\tilde{q_2},\tilde{q_3},\tilde{q_4}$ given by 
$$\tilde{f_1}=C_1,~\tilde{f_2}=C_2,~\tilde{f_3}=C_3,~\tilde{f_4}=C_4.$$

Here $(C_1,C_2,C_3,C_4)$ is a generic point of $A=\mathbb{G}_m^4$. 

Note that equations for the curve $E$ can be written as
\begin{equation} \label{E}
\begin{split}
(1-q_1-q_2)(1-q_3-q_4)=y\cdot \Big(C_1x_1x_2x_3\Big),~~~~~~~~~~~~~~~~~~~~~~~~~~~~~~~~~~~\\
q_1q_3=(y-1)\cdot \Big(C_2(x_1-1)(x_2-1)(x_3-1)\Big),~~~~~~~~~~~~~~~~~~~~~~~~~~~~~~~~~\\
q_2q_4=(y-t)\cdot \Big(C_3(x_1-t)(x_2-t)(x_3-t)\Big),~~~~~~~~~~~~~~~~~~~~~~~~~~~~~~~~~~\\
\Big(q_1q_2+\frac{(x_1-t)(x_2-t)(y-t)}{t(t-1)^2}q_1+\frac{(x_1-1)(x_2-1)(y-1)}{(t-1)^2}q_2\Big) \cdot~~~~~~~~~~~~~~~ \\\Big(q_3q_4+\frac{(y-t)(x_3-t)(z-t)}{t(t-1)^2}q_3+\frac{(y-1)(x_3-1)(z-1)}{(t-1)^2}q_4\Big)=~~~~~~~~~~~~~~~\\
y(y-1)(y-t)\cdot \Big(C_1C_2C_3C_4x_1x_2x_3(x_1-1)(x_2-1)(x_3-1)(x_1-t)(x_2-t)(x_3-t)\Big)
\end{split}
\end{equation}
and equations for the curve $\tilde{E}$ are obtained from the equations (\ref{E}) by interchanging $x_2$ and $x_3$, and replacing $y,q_1,q_2,q_3,q_4$ by $\tilde{y},\tilde{q_1},\tilde{q_2},\tilde{q_3},\tilde{q_4}$.
One can show by direct computation that $E$ and $\tilde{E}$ are elliptic curves with the same $j$-invariant. To show this we solve the first three equations in the system (\ref{E}) with respect to $q_1,q_2,y$ and substitute the result into the forth equation. We obtain a plane cubic curve with coordinates $q_3,q_4$, its genus and $j$-invariant can be computed in a usual way. After that we observe that $j$-invariant is symmetric with respect to $x_1,x_2,x_3$.

Observe that there exist rational points
$$(y,q_1,q_2,q_3,q_4)=(0,-\frac{C_2(t-1)}{z-1},\frac{C_3t(t-1)}{z-t},\frac{z-1}{t-1},-\frac{z-t}{t-1})\in E,$$
$$(\tilde{y},\tilde{q_1},\tilde{q_2},\tilde{q_3},\tilde{q_4})=(0,-\frac{C_2(t-1)}{z-1},\frac{C_3t(t-1)}{z-t},\frac{z-1}{t-1},-\frac{z-t}{t-1})\in \tilde{E}.$$
This gives a birational mapping $\mu:~E\to \tilde{E}$ such that $$\mu(0,-\frac{C_2(t-1)}{z-1},\frac{C_3t(t-1)}{z-t},\frac{z-1}{t-1},-\frac{z-t}{t-1})=(0,-\frac{C_2(t-1)}{z-1},\frac{C_3t(t-1)}{z-t},\frac{z-1}{t-1},-\frac{z-t}{t-1}).$$
One can also check that $\mu$ transforms $h~ dy\wedge dq_1\wedge dq_2\wedge dq_3\wedge dq_4$ to $\tilde{h}~ d\tilde{y}\wedge d\tilde{q_1}\wedge d\tilde{q_2}\wedge d\tilde{q_3}\wedge d\tilde{q_4}$.   $\square$

\subsection{A hypothetical example of a generalized product}

Here we suggest a hypothetical and more complicated example related to the generalized products described in Section 4.1. We still have $X=\A^1$ with auxiliary 
integration over $(n+1)$-dimensional cycle, but our commutative  associative operation has $n+1$ inputs and $n$ outputs.

{\bf Example 2.4.1.} Let 
\begin{equation}\label{Kn}
K(x_1,...,x_{2n+1},q_1,...,q_{n+1})=(x_1...~ x_{2n+1})^su_1^{-k_1}...~u_{n+1}^{-k_{n+1}}q_1^{2k_1-1}...~q_{n+1}^{2k_{n+1}-1}\times
\end{equation}
$$\Big(1+q_1+...+q_{n+1}\Big)^{s-k_1-...-k_{n+2}}\Big(1+\frac{u_1}{q_1}+...+\frac{u_{n+1}}{q_{n+1}}\Big)^{s+k_1+...+k_{n+2}}$$
where
$$u_1=\frac{(x_1-1)...(x_{2n+1}-1)t_1^2...~t_n^2}{x_1...~x_{2n+1}(t_1-1)^2...(t_n-1)^2},$$
$$u_{i+1}=\frac{(x_1-t_i)...(x_{2n+1}-t_i)\prod_{j\ne i} t_j^2}{x_1...~x_{2n+1}(t_i-1)^2\prod_{j\ne i}(t_j-t_i)^2},~i=1,...,n.$$
Here $t_1,...,t_n\ne 0,1$ are pairwise distinct parameters and $s,k_1,...,k_{n+2}$ are arbitrary constants. 

{\bf Conjecture 2.4.1.} There exists a birational mapping
$$(y_1,...,y_n,q_1,...,q_{2n+2})\to (\tilde{y}_1,...,\tilde{y}_n,\tilde{q}_1,...,\tilde{q}_{2n+2})$$
which transforms the l.h.s. of the equation
$$K(x_1,...,x_{n+1},y_1,...,y_n,q_1,...,q_{n+1})K(x_{n+2},y_1,...,y_n,z_1,...,z_n,q_{n+2},...,q_{2n+2})\times$$
$$dy_1\wedge ...\wedge dy_n\wedge dq_1\wedge...\wedge dq_{2n+2}=$$
$$\pm K(x_1,...,x_{n+2},\tilde{y}_1,...,\tilde{y}_n,\tilde{q}_1,...,\tilde{q}_{n+1})K(x_{n+1},\tilde{y}_1,...,\tilde{y}_n,z_1,...,z_n,\tilde{q}_{n+2},...,\tilde{q}_{2n+2})\times$$
$$d\tilde{y}_1\wedge ...\wedge d\tilde{y}_n\wedge d\tilde{q}_1\wedge...\wedge d\tilde{q}_{2n+2}=$$
to its r.h.s. Here the r.h.s. is obtained from the l.h.s. by interchanging $x_{n+1}$ and $x_{n+2}$, and replacing the variables 
$y_1,...,y_n,q_1,...,q_{2n+2}$ by $\tilde{y}_1,...,\tilde{y}_n,\tilde{q}_1,...,\tilde{q}_{2n+2}$.

Notice that in the Conjecture above $n>1$ because if $n=1$, then this family of kernels is essentially the same as in Example 5, so for $n=1$ this is proved.

The kernel in this example can be written in more symmetric form 
$$\int K(x_1,...,x_{2n+1},q_1,...,q_{n+1})dq_1...dq_{n+1}=$$ 
$$\int~~~ \prod_{i=1}^{n+2}(q_{i,+}^{+k_i}q_{i,-}^{-k_i})\cdot \Big(\sum_{i=1}^{n+2}q_{i,+}\Big)^{s-k_1-...-k_{n+2}}\Big(\sum_{i=1}^{n+2}q_{i,-}\Big)^{s+k_1+...+k_{n+2}}\cdot\frac{\prod_{i=1}^{n+2}\frac{dq_{i,+}}{q_{i,+}}}{\frac{d\lambda}{\lambda}}$$
where we integrate over $(n+1)$-dimensional torus $\mathbb{G}_m^{n+2}/ (\mathbb{G}_m)_{diag}$ with coordinates $q_{1,+},...,q_{n+2,+}$ factorized by the diagonal action 
$q_{i,+}\mapsto \lambda q_{i,+}$ of $\mathbb{G}_m$. Variables $q_{1,-},...,q_{n+2,-}$ are defined by 
$$q_{i,+}q_{i,-}=\frac{\prod_{\alpha=1}^{2n+1}(x_{\alpha}-t_i)}{P^{\prime}(t_i)^2},~~~i=1,...,n+2$$
where $P(u)=\prod_{i=1}^{n+2}(u-t_i)$.
Here we use arbitrary pairwise distinct parameters $t_1,...,t_{n+2}$, the previous formula for the kernel $K(x_1,...,x_{2n+1},q_1,...,q_{n+1})$ corresponds to the 
choice $t_{n+1}=0,~t_{n+2}=1$. 

The expression which we integrate is invariant under the diagonal action of $\mathbb{G}_m$. It is also invariant with respect to the group $S_{n+2}$ acting on 
variables $q_{1,+},...,q_{n+2,+}$ by permutations. One can achieve invariance with respect to a larger group $S_{n+2}\ltimes (\Z/2)^{n+2}$ (where the $i$th generator 
of $(\Z/2)^{n+2}$ acts by interchanging of $q_{i,+}$ and $q_{i,-}$) by adding one extra variable of integration (and loosing geometric convergence).
$$\int K(x_1,...,x_{2n+1},q_{1},...,q_{n+2})dq_{1}...dq_{n+2}=C\int~ \prod_{i=1}^{n+2}(q_{i,+}^{+k_i}q_{i,-}^{-k_i})\cdot \Big(\sum_{i=1}^{n+2}(q_{i,+}+q_{i,-})\Big)^{2s}\cdot\prod_{i=1}^{n+2}\frac{dq_{i,+}}{q_{i,+}}$$
where we integrate over the $(n+2)$-dimensional torus with coordinates $q_{1,+},...,q_{n+2,+}$. Here $C$ is independent of $x_1,...,x_{2n+1}$. 

Indeed, our two expressions for the kernel can be written as 
$$F_1=\int~ \Bigg(\sum_{i=1}^{n+2}a_iq_i+\sum_{i=1}^{n+2}\frac{b_i}{q_i}\Bigg)^{\alpha+\beta}\cdot \prod_{i=1}^{n+2}q_i^{\lambda_i}\cdot \prod_{i=1}^{n+2}\frac{dq_i}{q_i},$$
$$F_2=\int~ \Bigg(\sum_{i=1}^{n+2}a_iq_i\Bigg)^{\alpha}\cdot \Bigg(\sum_{i=1}^{n+2}\frac{b_i}{q_i}\Bigg)^{\beta}\cdot \prod_{i=1}^{n+2}q_i^{\lambda_i}\cdot \prod_{i=1}^{n+2}\frac{dq_i}{q_i}$$
for some $a_i,b_i,\lambda_i,\alpha,\beta$ such that 
$$\sum_{i=1}^{n+2}\lambda_i=\beta-\alpha.$$
Notice that the expressions $F_1,F_2$ both satisfy the same holonomic system of differential equations\footnote{Such expressions and the corresponding $D$-modules belong to 
a class of A-hypergeometric functions \cite{GG}.}:
$$\Big(a_i\frac{\partial}{\partial a_i}-b_i\frac{\partial}{\partial b_i}+\lambda_i\Big)F=0,~~~i=1,...,n+2,$$
\begin{equation}\label{Dn}
\Bigg(\sum_{i=1}^{n+2}\Big(a_i\frac{\partial}{\partial a_i}+b_i\frac{\partial}{\partial b_i}\Big)-\alpha-\beta\Bigg)F=0,
\end{equation}
$$\frac{\partial^2F}{\partial a_1\partial b_1}=\frac{\partial^2F}{\partial a_2\partial b_2}=...=\frac{\partial^2F}{\partial a_{n+2}\partial b_{n+2}}$$
which means that they should be in a sense equal.

The equality $F_2=CF_1$ can be also shown as follows.
Introduce an auxiliary function 
$$\phi(A,B)=\int \delta\Big(\sum_{i=1}^{n+2}a_iq_i-A\Big)\cdot \delta\Big(\sum_{i=1}^{n+2}\frac{b_i}{q_i}-B\Big)\cdot \prod_{i=1}^{n+2}q_i^{\lambda_i}\cdot \prod_{i=1}^{n+2}\frac{dq_i}{q_i}.$$
The function $\phi(A,B)$ is homogeneous 
$$\phi(qA,q^{-1}B)=q^{\beta-\alpha}\phi(A,B)~~~\text{for all}~~~q\ne 0,$$ 
so we can write $$\phi(A,B)=\phi_0(AB)\cdot A^{\beta-\alpha}$$ 
where $\phi_0(u)=\phi(1,u)$.
We have 
$$F_1=\int (A+B)^{\alpha+\beta}\phi(A,B)dAdB=\int (A+B)^{\alpha+\beta}\phi_0(AB)\cdot A^{\beta-\alpha}dAdB=$$ 
$$\int \Big(A+\frac{u}{A}\Big)^{\alpha+\beta}\phi_0(u)A^{\beta-\alpha}du\frac{dA}{A}=\int u^{\beta}\phi_0(u)du\cdot \int \frac{(1+v^2)^{\alpha+\beta}}{2v^{2\alpha}}\frac{dv}{v}$$
where we made substitutions $B=\frac{u}{A}$ and $A=v\sqrt{u}$.

Similar computation for $F_2$ gives
$$F_2=\int A^{\alpha}B^{\beta}\phi(A,B)dAdB=\int A^{\alpha}\Big(\frac{u}{A}\Big)^{\beta}\phi_0(u)A^{\beta-\alpha}du\frac{dA}{A}=\int u^{\beta}\phi_0(u)du\cdot \int \frac{dA}{A}.$$
Therefore, $F_1$ and $F_2$ are both equal to the integral $\int u^{\beta}\phi_0(u)du$ multiplied by a constant independent of $a_i,b_i,~i=1,...,n+2$. Moreover, 
removing the integral $\int \frac{dA}{A}$ from our final expression for $F_2$ we cure its geometrical divergence. 

\subsection{Some examples with 3 inputs and one output}

Here we give examples of products $\mu_3$ with 3 inputs and one output with look a bit pathological but still satisfy an analog of  associativity and commutativity conditions:  $\mu_3$ is $S_3$-invariant and $\mu_3\circ (\mu_3\otimes id)$ is $S_5$-invariant.

{\bf Example 2.5.1.} Let 
$$K(x_1,x_2,x_3,y)=e^{c(x_1x_2x_3y+y)}~\frac{1}{y}$$
where $c$ is a constant. Then we have
$$K(x_1,x_2,x_3,y)~K(y,x_4,x_5,z)~dy=K(x_1,x_2,x_4,\tilde{y})~K(\tilde{y},x_3,x_5,z)~d\tilde{y}$$
if $\tilde{y}=\frac{x_1x_2x_3+x_4x_5z+1}{x_1x_2x_4+x_3x_5z+1}y.$

This example is similar to Example 1 and it is degenerate in the same sense. 

{\bf Example 2.5.2.} Let 
$$K(x_1,x_2,x_3,y)=\Big(\frac{1}{y}+\frac{1}{y}\sqrt{1+x_1x_2x_3y}\Big)^c~\frac{1}{\sqrt{1+x_1x_2x_3y}}$$
where $c$ is an arbitrary constant. Then 
\begin{equation} \label{ass4}
K(x_1,x_2,x_3,y)~K(y,x_4,x_5,z)~dy=K(x_1,x_2,x_4,\tilde{y})~K(\tilde{y},x_3,x_5,z)~d\tilde{y}
\end{equation}
where $\tilde{y}$ is a function in $x_1,x_2,x_3,x_4,x_5,y,z$ independent of $c$. 

More precisely, let $C_1\subset \A^3$ be a rational affine curve given by\footnote{Here $w_{123},w_{45},y$ are affine coordinates on $\A^3$ and $x_1,x_2,x_3,x_4,x_5,z$ are parameters.} 
$$w_{123}^2=1+x_1x_2x_3y,~~~w_{45}^2=1+yx_4x_5z$$
and $C_2\subset \A^3$ be a rational affine curve given by\footnote{Here $w_{124},w_{35},\tilde{y}$ are affine coordinates on $\A^3$ and $x_1,x_2,x_3,x_4,x_5,z$ are parameters.} 
$$w_{124}^2=1+x_1x_2x_4\tilde{y},~~~w_{35}^2=1+\tilde{y}x_3x_5z.$$
Then the formulas
$$w_{124}=\frac{(x_1x_2-x_5z)x_4w_{123}+(x_3-x_4)x_1x_2w_{45}}{x_1x_2x_3-x_4x_5z},$$
$$w_{35}=\frac{(x_3-x_4)x_5zw_{123}+(x_1x_2-x_5z)x_3w_{45}}{x_1x_2x_3-x_4x_5z},$$
$$\tilde{y}=\frac{2(x_3-x_4)(x_1x_2-x_5z)(w_{123}w_{45}-1)+x_3x_4(x_1x_2-x_5z)^2y+x_1x_2x_5z(x_3-x_4)^2y}{(x_1x_2x_3-x_4x_5z)^2}$$
define a birational mapping $C_1\to C_2$ which transforms the l.h.s. of (\ref{ass4}) to its r.h.s.

\subsection{On $S_n$-covariance of higher compositions of kernels}

If $K$ is a multiplication kernel of birational type on variety $X/k$, then for any $n\geq 2$ and any planar binary rooted tree $T$ with $n$ leaves (or, equivalently, 
choice of bracketing on a product of $n$ symbols) we get a variety $Z_T$ which maps to $X^{n+1}\times A$ and is endowed with a volume element along fibers for the 
projection to the last factor $X$. This variety 
$Z_T$ corresponds to the kernel of $n$-fold product with the chosen bracketing. Associativity implies that there exists a birational identification of varieties $Z_T$ for 
different trees $T$. Commutativity of our kernel implies that one can lift each permutation $\sigma\in S_n$ of first $n$ factors in $X^{n+1}$, to the birational 
identification of varieties $Z_T$. These identifications, however, are not compatible in general. It is desirable to lift these identifications to a $S_n$-action and 
construct just one variety $Z_n$ (with an action of $S_n$) which maps to $X^{n+1}\times A$ in a $S_n$-covariant way. 

One of the reasons why it is desirable is the following. Let $k$ be a local field, and assume that the kernel is geometrically convergent. Then we expect that an appropriate space 
${\mathcal F}(X(k))$ of complex-valued functions on $X(k)$ to be a commutative associative algebra with product $*$ given by our kernel $K$. For any finite extension $k^{\prime}\supset k$ we get another algebra $({\mathcal F}(X(k^{\prime})),*)$. The existence of $S_n$-action on $Z_n$ for $n=deg(k^{\prime}/k)$ gives rise to a 
homomorphism 
\begin{equation}\label{ftof}
({\mathcal F}(X(k^{\prime})),*)\to ({\mathcal F}(X(k)),*).
\end{equation}
Namely, extension $k^{\prime}\supset k$ gives a transitive action of $Gal(\bar{k}/k)$ on an $n$-elements set (the set of embeddings $k^{\prime}\hookrightarrow  \bar{k}$ 
over $k$), hence a homomorphism $\rho: Gal(\bar{k}/k)\to S_n$, up to conjugation. Using $\rho$ we can define a twisted form $Z_{n,\rho}$ of $Z_n$. Recall that the set 
$Z_{n,\rho}(\bar{k})$ of $\bar{k}$-points of $Z_{n,\rho}$ coincides with $Z_n(\bar{k})$, but the action of $Gal(\bar{k}/k)$ is twisted by $\rho$. 
Variety $Z_{n,\rho}$ maps to the product $(X^n)_{\rho}\times X\times A$ where $(X^n)_{\rho}$ is the twisted forms of 
$X^n$, i.e. the Weil restriction $Res_{k^{\prime}/k}(X)$. Recall that the latter is a variety over $k$ such that its set of $k$-points is $X(k^{\prime})$. The twisted 
multiplication kernel gives a map (\ref{ftof}). One can check that it is a homomorphism of algebras. 

In the case of varieties over finite fields (see setup 4 in Section 1.2) there are homomorphisms (\ref{ftof}) for extensions $k^{\prime}=\mathbb{F}_{q^m}\supset 
k=\mathbb{F}_q$ of finite fields, even without the action of $S_n$ on $Z_n$. It looks plausible that in this case there is still an action of $S_n$ on the corresponding 
motivic constructible sheaves. The tower of homomorphisms (\ref{ftof}), or dually, maps (which are inclusions for finite fields)
$$Spec~({\mathcal F}(X(k)),*)\hookrightarrow Spec~({\mathcal F}(X(k^{\prime})),*)$$
played an essential role in \cite{K}. The inductive limit 
$$\lim_{\underset{m}{\longrightarrow}} Spec~({\mathcal F}(X(\mathbb{F}_{q^m}),*)(\overline{\Q})$$
is an infinite countable set endowed with  
commuting actions of two Galois groups: $Gal(\overline{\mathbb{F}_q}/\mathbb {F}_q)$ and $Gal(\overline{\Q}^{CM}/\Q)$.

Unfortunately, in general it seems to be impossible to lift $S_n$ action from $X^n\times X\times A$ to $Z_n$. Let us describe this problem in the case of Example 2.2.2. 
Every binary rooted tree $T$ with $n\geq 3$ leaves gives a family of Calabi-Yau varieties $V^T_{x_1,...,x_{n+1},t}$ of dimension $n-3$ depending on $n+2$ parameters 
$(x_1,...,x_{n+1},t)\in \big(\mathbb{G}_m\big)^{n+1}\times\mathbb{G}_a$ in the following way. Variety $V^T_{x_1,...,x_{n+1},t}$ is a hypersurface in a toric variety, and it is given by equation in variables 
$(y_1,...,y_{n-2})\in\mathbb{G}_m^{n-2}$ 
$$\sum_{v\in \{\text{vertices of } T\}}f_v=t$$
where parameters $x_1,...,x_n$ are attached to leaves of $T$, parameter $x_{n+1}$ is attached to the root of $T$, and parameters $y_1,...y_{n-2}$ are attached to inner 
edges of $T$. For each vertex $v$ we define 
$$f_v=f(z_1,z_2,z_3)=z_1z_2z_3+\frac{z_1}{z_2z_3}+\frac{z_2}{z_1z_3}+\frac{z_3}{z_1z_2}$$
where $z_1,z_2,z_3$ are variables attached to three edges adjacent to $v$. For example, if 

$
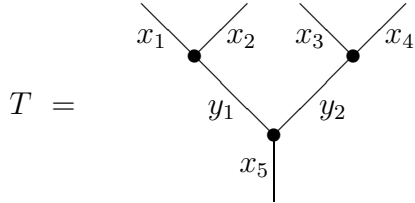
\begin{picture}(150,80)(-100,0)
 \put(100,20){\circle*{5}}
 \put(100,20){\line(0,-1){25}}
 \put(100,20){\line(-1,1){50}}
 \put(100,20){\line(1,1){50}}
 \put(130,50){\circle*{5}}
 \put(70,50){\circle*{5}}
 \put(130,50){\line(-1,1){20}} \put(70,50){\line(1,1){20}}
 \put(48,55){$x_1$}
 \put(108,55){$x_3$}
 \put(82,55){$x_2$}
 \put(142,55){$x_4$}
 \put(75,28){$y_1$}
 \put(117,28){$y_2$}
 \put(87,6){$x_5$}
 \put(0,28){$T\,\,=$}
 \end{picture}$
 \\
 \\
 then $V^T_{x_1,...,x_5,t}$ is an elliptic curve given by 
 $$f(x_1,x_2,y_1)+f(x_3,x_4,y_2)+f(y_1,y_2,x_5)=t.$$
 Variable $t$ is the coordinate on $A=\mathbb{G}_a=\A^1$. 
 
 The volume element (up to a sign) on $V^T_{x_1,...,x_{n+1},t}$ is given by 
 $$\frac{\wedge_{i=1}^{n-2}d\log y_i}{d(\sum_v f_v)}.$$
 The integral operator corresponding to $Z_T$ is given by a density on $X^{n+1}=(\mathbb{G}_m)^{n+1}$ with coordinates 
 $x_1,...,x_{n+1}$. This density is the Fourier transform of function $t\mapsto vol(V^T_{x_1,...,x_{n+1},t})$ in variable $t$. One can construct a birational 
 identification $V^T_{x_1,...,x_{n+1},t}\sim V^{T^{\prime}}_{x_1,...,x_{n+1},t}$ there $T^{\prime}$ is obtained by a flip from $T$:
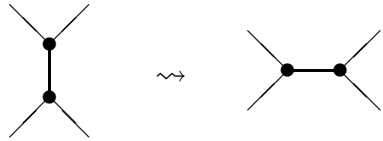
 $$\begin{picture}(150,60)(0,0)
 \put(5,5){\line(1,1){15}}
 \put(20,20){\circle*{5}}
 \put(20,20){\line(1,-1){15}}
 \put(20,20){\line(0,1){20}}
 \put(20,40){\circle*{5}}
 \put(20,40){\line(-1,1){15}}
 \put(20,40){\line(1,1){15}}
 \put(60,25){$\rightsquigarrow$}
 \put(110,30){\circle*{5}}
 \put(110,30){\line(-1,1){15}}
 \put(110,30){\line(-1,-1){15}}
 \put(110,30){\line(1,0){20}}
 \put(130,30){\circle*{5}}
  \put(130,30){\line(1,1){15}}
 \put(130,30){\line(1,-1){15}}
 \end{picture}$$
 The composition of 5 flips corresponding to pentagon relation is a non-trivial automorphism of K3 surface $V^T_{x_1,...,x_5,t}$. 
 
 It might be possible to cure this problem by increasing the dimension of $Z_n$ (or equivalently, adding auxiliary variables for the integration). This might require 
 an extension of our formalism. One such possibility is discussed in the next Section. 

\subsection{Direct images of cyclic {\it D-}modules at a generic point}

In this section we propose a mixed birational/$D$-module type formalism using which one can speak about multiplication kernels. In this formalism, like in Section 2.1, 
one considers algebraic varieties only up to birational equivalence. On the other hand, we deal here with cyclic $D$-modules. Roughly speaking, we encode a ``function'' 
$f$ on algebraic variety $X$ by a system of algebraic linear differential equations satisfied by $f$ at the generic point of $X$. 
The most interesting operation with these objects is the direct image (which we denote by $\pi_*$) defined below, which informally corresponds to the integration over 
an unspecified cycle. Almost all examples in our paper can be understood in this mixed formalism. 

Let $X$ be a smooth algebraic variety over field $k$ of characteristic zero, endowed with a line bundle ${\mathcal L}$. Denote by $\text{Diff}_{{\mathcal L},rat}$ the algebra of 
differential operators in ${\mathcal L}$ with coefficients in rational functions on $X$. Denote by $\I_{\mathcal L}$ the set of left ideals in 
$\text{Diff}_{{\mathcal L},rat}$. We can think about elements of $\I_{\mathcal L}$ as systems of differential equations on a ``section'' of ${\mathcal L}$ (e.g. 
analytic germ of a section for $k=\C$). For any $I\in \I_{\mathcal L}$ the quotient $\text{Diff}_{{\mathcal L},rat}/I$ is a cyclic $\text{Diff}_{{\mathcal L},rat}$-module.
 In this paper we deal only with examples where this cyclic module is holonomic. It is clear that both $\text{Diff}_{{\mathcal L},rat}$ and ${\mathcal L}$ depends 
 only on the birational type of $X$ (and also on ${\mathcal L}$), i.e. only on the field $k(X)$ of rational functions, together with a 1-dimensional module over it.
 
 There are three basic operations:
 
 {\bf 1)} For a dominant map $\pi:Y\to X$, and ${\mathcal L}$ on X, we have a pullback $\pi^*:\I_{\mathcal L}\to\I_{\pi^*{\mathcal L}}$.
 
 {\bf 2)} In the same situation, we have pushforward $\pi_*:\I_{~\pi^*{\mathcal L}\otimes K_{Y/X}}\to\I_{\mathcal L}$.
 
 {\bf 3)} For a collection of line bundles ${\mathcal L}_1,...,{\mathcal L}_n$ on $X$ we have product $\prod_i \I_{{\mathcal L}_i}\to \I_{\otimes_i {\mathcal L}_i}$.
 
 In order to explain {\bf 1)} and {\bf 2)} it is convenient to have a non-linear flat connection on fibration $\pi: Y\to X$. We claim that there exists a Zariski open 
 dense set $Y^{\prime}\subset Y$ such that $\pi|_{Y^{\prime}}$ is a submersion, and a non-linear flat connection $\nabla$ on $\pi|_{Y^{\prime}}:Y^{\prime}\to X$, i.e. an integrable 
 distribution on $Y^{\prime}$ transversal to fibers of $\pi$. 
 
 Indeed, passing to an open dense subset of $Y$ we may assume that fibers of $\pi$ are smooth locally closed subvarieties of $\A^N$ for some $N$. Let us choose a generic 
 affine projection $\A^N\to\A^n$ where $n=\dim Y-\dim X$. Then the fibers $\pi^{-1}(x)$ for $x\in X$ are ramified coverings over the {\it same} affine space $\A^n$, 
 hence are locally identified in etale topology. This gives a non-linear flat connection outside of the ramification loci $\pi^{-1}(x)\to\A^n,~x\in X$. Any non-linear 
 flat connection $\nabla$ gives rise to map of algebras $\nabla^*: \text{Diff}_{{\mathcal L},rat} \to \text{Diff}_{\pi^*{\mathcal L},rat}$. After making the choice of 
 connection $\nabla$, the maps $\pi^*,\pi_*$ in {\bf 1)}, {\bf 2)} are defined as follows. In the case {\bf 1)}, it is convenient to assume ${\mathcal L}={\mathcal O}_X$ 
 (this can be achieved by passing to an open dense set $X^{\prime}\subset X$). For an ideal $I\in \I_{{\mathcal O}_X}$ we define $\pi^*I$ to be the left ideal in 
 $\I_{{\mathcal O}_X}$ generated by elements $\nabla^*(P)$ for $P\in\I$, and by vector fields along fiber of $\pi$.
 
 In the case {\bf 2)} it is convenient to assume ${\mathcal L}=K_X$, and the left ideals in $\text{Diff}_{{\mathcal L},rat}$ are the same as right ideals in 
 $\text{Diff}_{{\mathcal O}_X,rat}=\text{Diff}_{{\mathcal L},rat}^{op}$. Let $I$ be a right ideal in $\text{Diff}_{Y,rat}$ which can be though of as an annihilator of a 
 ``volume form'' on $Y$ (e.g. analytic volume form in some domain for $k=\C$). We want to define $\pi_*I$ to be the right ideal in $\text{Diff}_{X,rat}$ consisting of 
 differential operators annihilating formal integrals $\pi_*(vol_Y)$. By definition, $P\in \text{Diff}_{X,rat}$ belongs to $\pi_*I$ iff $\nabla^*(P)~ \text{mod} ~\I_Y\in 
 \I_Y\backslash\text{Diff}_{Y,rat}$ belongs to the image of the right submodule in $\text{Diff}_{Y,rat}$ generated by vector fields along fibers of $\pi$. 
 
 In both cases {\bf 1)}, {\bf 2)} the result of operations does not depend on the choice of flat connection $\nabla$.
 
 In the case {\bf 3)}, we assume that all ${\mathcal L}_1,...,{\mathcal L}_n$ are trivialized. It is enough to consider the case $n=2$. Define a linear map 
 $\Delta: \text{Diff}_{{\mathcal O}_X,rat}\to\text{Diff}_{{\mathcal O}_X,rat}\otimes_{k(X)}\text{Diff}_{{\mathcal O}_X,rat}$ (here both factors on the right are understood
  as left $k(X)$-modules), by the condition $\Delta P=\sum_{\alpha}P^{(1)}_{\alpha}\otimes P^{(2)}_{\alpha}$ if for any $f,g\in k(X)$ we have 
  $P(fg)=\sum_{\alpha}P^{(1)}_{\alpha}f\cdot P^{(2)}_{\alpha}g$. For two left ideals $I_1,I_2\in\I_{{\mathcal O}_X}$ we define their ``product'' $J$ by 
  $J=\{P\in \text{Diff}_{{\mathcal O}_X,rat};~ \Delta P\in I_1\otimes \text{Diff}_{{\mathcal O}_X,rat}+\text{Diff}_{{\mathcal O}_X,rat}\otimes I_2\}$. 
  This is the ideal annihilating product $f_1f_2$ where $f_1,f_2$ are general analytic functions annihilated by $I_1$, $I_2$ respectively. 
  
\section{Semi-classical kernels and their quantization}

\subsection{Semi-classical kernels associated with classical integrable systems}

{\bf Definition 3.1.1.} Let $(M,\omega)$ be a symplectic manifold. A semi-classical multiplication kernel is a tuple $((M,\omega),N,P)$ 
where $N\subset (M,\omega)$, $P\subset (M,-\omega)\times (M,-\omega)\times (M,\omega)$ are Lagrangian submanifolds with the following properties:

{\bf 1.} Projections $\pi_i:P\to M,~i=1,2,3$ are submersions.

{\bf 2.} $P$ is symmetric with respect to the interchanging of the first two components in $M^3$.

{\bf 3.} The correspondence $M\times M\to M$ given by $P$ is associative. 

{\bf 4.} The map $(\pi_2\times\pi_3)|_{\pi_1^{-1}(N)\cap P}$ is an inclusion whose image is open dense in the diagonal Lagrangian submanifold $M_{diag}=\{(m,m);~m\in M\}\subset (M,-\omega)\times (M,\omega)$.

{\bf Remark 3.1.1.} In other words, $(M,\omega)$ is a commutative monoid in  the  symplectic monoidal category introduced by A. Weinstein.  The product in this commutative 
monoid is defined by $P$ and the neutral element is defined by $N$. In general, the objects of this category are symplectic manifolds (in algebraic or $C^{\infty}$ sense) 
and 
morphisms $Hom((M_1,\omega_1),(M_2,\omega_2))$ are defined as Lagrangian submanifolds in $(M_1,-\omega_1)\times (M_2,\omega_2)$. The tensor product is given by the 
usual product of symplectic manifolds. Notice that there are transversality problems in the definition of the composition as the composition of correspondences. 
Strictly speaking, the above definition is a rough sketch, a first approximation to a not yet found more satisfactory and rigorous notion. 

{\bf Remark 3.1.2.} The constraint on Lagrangian subvariety $P$ to be smooth (i.e. to be a manifold) seems to be not completely natural. The closure of $P$ in $M^3$ is usually 
singular. 

It seems plausible that semi-classical multiplication kernels are essentially the same as Liouville integrable systems endowed with a Lagrangian section. 

Let $(M,\omega)$ be a symplectic manifold with a structure of Liouville integrable system. In other words, we have a fibration  
$$f: M\to B$$ 
where $\dim B=\frac{1}{2} \dim M$ and such that pullback of functions on $B$ Poisson commute. Then $L_b=f^{-1}(b)\subset M$ is a 
Lagrangian submanifold of $M$ for generic point $b\in B$, endowed with a natural affine structure (a torsion-free flat connection on the tangent bundle $TL_b$). 
Let us assume for simplicity that the generic fiber $L_b$ is connected and compact, and moreover a Lagrangian section $\sigma:B\to M$, $f\circ\sigma=id$ is chosen. 
Then $L_b$ will have a structure of a commutative group with the identity element given by $\sigma(b)\in L_b$. In fact, $L_b$ is an abelian variety\footnote{In practice,
fibers of $f$ are often noncompact and admit compactification to abelian varieties. On the other hand, in the case $\dim M=2$, there are many integrable systems 
for which the generic fiber of $f$ is a punctured curve of genus $g>1$, e.g. one can take $M=T^*\A^1=\A^2$ and $B=\A^1$, with the map $f$ given by a polynomial in 2 
variables of sufficiently high degree.}. Define $P_b\subset L_b\times L_b\times L_b\subset M\times M\times M$ 
by $P_b=\{(u,v,u+v);u,v\in L_b=f^{-1}(b)\}$. Then one can see that $P=\cup_{b\in B}P_b$ is a Lagrangian submanifold in $(M,-\omega)\times (M,-\omega)\times (M,\omega).$ 
In this way we obtain a semi-classical multiplication kernel given by the correspondence $P$, with the neutral element given by $N=\sigma(B)\subset M$. 

{\bf Remark 3.1.3.} If we do not choose section $\sigma$ (or submanifold $N$), then we will have naturally only a structure of an abelian torsor on $L_b,b\in B$. 
Instead of multiplication $P\subset (M,-\omega)\times (M,-\omega)\times (M,\omega)$ we will have a manifold $\{(u_1,u_2,u_3,u_4)\in M^4; f(u_1)=f(u_2)=f(u_3)=f(u_4), u_1+u_2=u_3+u_4\}$ which is a Lagrangian submanifold in $(M,-\omega)\times (M,-\omega)\times (M,\omega)\times (M,\omega).$ This is a structure similar to one in 
Example 2, Section 1.5 where we consider commuting operators with simple joint spectrum and without the choice of normalization for the eigenfunctions.

Let $M=T^{*}X$ be cotangent bundle with the canonical symplectic structure where $X$ is an algebraic variety. Furthermore, to simplify notations we assume that 
$X=\A^1$ with an affine coordinate $x$. In this case $M\cong \A^2$ with canonical coordinates $x,p$ and $\omega=dx\wedge dp$.
Let $f:~M\to B$ is given by the formula $(x,p)\mapsto f(x,p)$ where $f(x,p)$ is a meromorphic function. In this case the graph of our commutative 
associative multiplication is defined by the system of equations of the form
\begin{equation} \label{clas}
f(x_1,p_1)=f(x_2,p_2)=f(x_3,p_3),~~~g(x_1,p_1,x_2,p_2,x_3,p_3)=0
\end{equation}
where $(x_1,p_1,x_2,p_2,x_3,p_3)$ are coordinates on $M\times M\times M$, and the first two equations mean that $f(m_1)=f(m_2)=f(m_3)=b$ for some $b\in B$. Here 
$(m_1,m_2,m_3)\in M\times M\times M$ and $(x_i,p_i)$ are coordinates of point $m_i$, $i=1,2,3$.

Notice that commutativity of our semi-classical multiplication kernel means that\footnote{Strictly speaking, we need only the equivalence $g(x_1,p_1,x_2,p_2,x_3,p_3)=0 \iff  g(x_2,p_2,x_1,p_1,x_3,p_3)=0$.} 
$$g(x_1,p_1,x_2,p_2,x_3,p_3)=g(x_2,p_2,x_1,p_1,x_3,p_3)$$
and associativity means that if we take the system of equations 
$$f(x_1,p_1)=f(x_2,p_2)=f(x_3,p_3)=f(x_4,p_4)=f(x_5,p_5),$$
$$g(x_1,p_1,x_2,p_2,x_3,p_3)=0,~g(x_3,p_3,x_4,p_4,x_5,p_5)=0$$
and exclude $x_3,p_3$ from it, the resulting affine variety will be symmetric with respect to interchanging of pairs $(x_2,p_2)$ and $(x_4,p_4)$.

Recall also that $P\subset (M,-\omega)\times (M,-\omega)\times (M,\omega)$ is a Lagrangian submanifold. This means that if $I$ is the ideal generated by the equations (\ref{clas}), then $$\{I,I\}\subset I$$ where $\{,\}$ are Poisson brackets defined by $\{p_1,x_1\}=\{p_2,x_2\}=-1,\{p_3,x_3\}=1$ and other brackets between coordinates are equal to zero. 

More generally, let $M=T^*\A^n$ with canonical coordinates $x_1,...,x_n,p_1,...,p_n$ be a phase space of an integrable system. Let $f_1(x_1,...,p_n),...,f_n(x_1,...,p_n)$ be commuting 
integrals of this integrable system. In this case Lagrangian submanifold\footnote{Here as usual $\omega=dx_1\wedge dp_1+...+dx_n\wedge dp_n$.} $P\subset (M,-\omega)\times (M,-\omega)\times (M,\omega)$ is defined by a system of equations in the variables $x_{1,i},...,x_{n,i},p_{1,i},...,p_{n,i},~i=1,2,3$ of the form
\begin{equation} \label{clasg}
\begin{split}
f_i(x_{1,1},...,x_{n,1},p_{1,1},...,p_{n,1})=f_i(x_{1,2},...,x_{n,2},p_{1,2},...,p_{n,2}),\\
f_i(x_{1,2},...,x_{n,2},p_{1,2},...,p_{n,2})=f_i(x_{1,3},...,x_{n,3},p_{1,3},...,p_{n,3}),\\
g_i(x_{1,1},...p_{n,3})=0,~i=1,...,n~~~~~~~~~~~~~~~~~
\end{split}
\end{equation}
for some functions $g_1,...,g_n$. This system should satisfy the following properties:

{\bf 1.} Let $I^{123}$ be the ideal generated by the equations (\ref{clasg}). Then $\{I^{123},I^{123}\}\subset I^{123}$ where $\{,\}$ is the canonical Poisson structure on the symplectic manifold $(M,-\omega)\times (M,-\omega)\times (M,\omega)$. This property means that $P$ is a Lagrangian submanifold in $(M,-\omega)\times (M,-\omega)\times (M,\omega)$.

{\bf 2.} Let $I^{213}$ be the ideal obtained from $I^{123}$ by interchanging of variables $x_{i,1},p_{i,1}$ and $x_{i,2},p_{i,2},~i=1,...,n.$ Then $I^{213}=I^{123}$. This property means commutativity of our semi-classical kernel.

{\bf 3.} Let $I^{345}$ be the ideal obtained from $I^{123}$ by replacing of variables $x_{i,1},p_{i,1},x_{i,2},p_{i,2},x_{i,3},p_{i,3}$ by the variables $x_{i,3},p_{i,3},x_{i,4},p_{i,4},x_{i,5},p_{i,5}$ for $i=1,...,n.$ Let $I^{12345}$ be the ideal generated by $I^{123}$ and $I^{345}$. 
Then the ideal obtained from $I^{12345}$ by excluding the variables $x_{i,3},p_{i,3},~i=1,...,n$ should be symmetric with respect to interchanging of variables $x_{i,2},p_{i,2}$ and $x_{i,4},p_{i,4},~i=1,...,n.$ This property means associativity of our semi-classical kernel.

\subsection{Quantization of semi-classical kernels and quantum integrable systems}

Let us discuss quantization of the picture above. Let us start with the case $M=T^*\A^1$. As usual, we replace the Poisson algebra in the variables $x_1,x_2,x_3,p_1,p_2,p_3$ by the ring 
of differential operators\footnote{Informally, we replace $p_1\mapsto \frac{d}{dx_1}$, $p_2\mapsto \frac{d}{dx_2}$, $p_3\mapsto -\frac{d}{dx_3}$ and choose a ``correct'' ordering of operators. Here we set Planck constant to one. The sign difference between derivatives by $x_1,x_2$ and $x_3$ reflects the sign difference in $\omega$ in the semi-classical case.} in $x_1,x_1,x_3,\frac{d}{dx_1},\frac{d}{dx_2},\frac{d}{dx_3}$. The system (\ref{clas}) should be replaced by a system of differential equations for a quantum multiplication kernel $K(x_1,x_2,x_3)$
\begin{equation} \label{quan}
\begin{split}
\hat{f}\Big(x_1,\frac{d}{dx_1}\Big)^*K(x_1,x_2,x_3)=\hat{f}\Big(x_2,\frac{d}{dx_2}\Big)^*K(x_1,x_2,x_3)=\hat{f}\Big(x_3,\frac{d}{dx_3}\Big)K(x_1,x_2,x_3),\\
\hat{g}\Big(x_1,\frac{d}{dx_1},x_2,\frac{d}{dx_2},x_3,\frac{d}{dx_3}\Big)K(x_1,x_2,x_3)=0~~~~~~~~~~~~~~~~~~~~~
\end{split}
\end{equation}
Here $\hat{f},\hat{g}$ are differential operators which are quantizations of $f,g$ and $^*$ is the anti-involution of the algebra of differential operators defined by $x_i^*=x_i,\Big(\frac{d}{dx_i}\Big)^*=-\frac{d}{dx_i}$, $(AB)^*=B^*A^*$. 

We assume that the left ideal $\hat{I}$ in the algebra of differential operators in $x_1,x_1,x_3,\frac{d}{dx_1},\frac{d}{dx_2},\frac{d}{dx_3}$ generated by the
equations (\ref{quan}), satisfies the property
$$[\hat{I},\hat{I}]\subset \hat{I}.$$
This is a quantization of the property $\{I,I\}\subset I$ in the semi-classical case. 

We also assume that $K(x_1,x_2,x_3)=K(x_2,x_1,x_3)$ and the 
differential operator 
$$\hat{g}\Big(x_1,\frac{d}{dx_1},x_2,\frac{d}{dx_2},x_3,\frac{d}{dx_3}\Big)$$ is symmetric with respect to interchanging of $x_1,x_2$.
This means that our quantum multiplication is commutative. 

It is less clear however how to lift associativity to the quantum case. We see the following two possibility:

{\bf P1.} There exists a non-zero solution $K(x_1,x_2,x_3)$ of the system (\ref{quan}) such that $K$ is a multiplication kernel of birational type.

{\bf P2.} Denote by $\hat{I}^{(123)}$ the left ideal in the algebra of differential operators generated by equations (\ref{quan}). Let 
$\hat{I}^{(345)}$ be obtained from $\hat{I}^{(123)}$ by replacing $x_1,x_2,x_3$ by $x_3,x_4,x_5$. Let 
$$\hat{I}^{(1245)}=(\hat{I}^{(123)}+\hat{I}^{(345)}+J)\cap R_{1245}$$ 
where $J$ is the {\it right} ideal in the ring of differential operators in variables $x_1,...,x_5$ generated by $\frac{d}{dx_3}$ 
and $R_{1245}$ is the ring of differential operators in $x_1,x_2,x_4,x_5$. We assume that $\hat{I}^{(1245)}$ is symmetric with respect to interchanging of $x_2$ and $x_4$. Notice that $\hat{I}^{(1245)}$ is the left ideal in the ring $R_{1245}$.

It is not clear how these two properties are related in general. Informally, if {\bf P1} holds, then 
\begin{equation} \label{c}
\int_{\Gamma}K(x_1,x_2,x_3)K(x_3,x_4,x_5)dx_3=\int_{\Gamma}K(x_1,x_4,x_3)K(x_3,x_2,x_5)dx_3
\end{equation}
for any cycle $\Gamma$. This means that the l.h.s. and the r.h.s. of (\ref{c}) should satisfy the same system of differential equations as a function in $x_1,x_2,x_4,x_5$. Therefore, {\bf P2} looks feasible but it is not clear how to make these considerations rigorous.

On the other hand, if {\bf P2} holds, then the l.h.s. and the r.h.s of (\ref{c}) should satisfy the same system of differential equations. In general this does not mean that the r.h.s. is obtained from the l.h.s. by a birational mapping, but it would be natural to expect this in examples.

The advantage of {\bf P2} is that it can in principle be verified algorithmically. In this paper, however, we concentrate on {\bf P1}.

Let us discuss quantization of semi-classilal kernels in a more general case. Consider a quantum integrable system defined by commuting differential operators $D_1,...,D_n$ in the variables $x_1,...,x_n$. A quantization of the system (\ref{clasg}) has a form
\begin{equation} \label{quang}
\begin{split}
D_{i,1}^*K=D_{i,2}^*K=D_{i,3}K,~i=1,...,n,\\
G_iK=0,~i=1,...,n~~~~~~~~~~
\end{split}
\end{equation}
for some differential operators $G_i$ in the variables $x_{i,j}$, $i=1,...,n,j=1,2,3$. Here $D_{i,j},~i=1,...,n,j=1,2,3$ are obtained from $D_i$ by replacing the variable $x_1,...,x_n$ by $x_{1,j},...,x_{n,j}$, and $K$ is a function in variables $x_{i,j}$. This system should satisfy the following properties:

{\bf 1.} Let $\hat{I}^{123}$ be the left ideal generated by the equations (\ref{quang}). Then $[\hat{I}^{123},\hat{I}^{123}]\subset \hat{I}^{123}$ where $[A,B]=AB-BA$ is the usual commutator of differential operators. This property means that the system (\ref{quang}) is holonomic.

{\bf 2.} Let $\hat{I}^{213}$ be the left ideal obtained from $\hat{I}^{123}$ by interchanging of variables $x_{i,1}$ and $x_{i,2},~i=1,...,n.$ Then $\hat{I}^{213}=\hat{I}^{123}$. This property means commutativity of our kernel.

{\bf 3.} Let $\hat{I}^{345}$ be the left ideal obtained from $\hat{I}^{123}$ by replacing of variables $x_{i,1},x_{i,2},x_{i,3}$ by the variables $x_{i,3},x_{i,4},x_{i,5}$ for $i=1,...,n.$ Let 
$$\hat{I}^{1245}=(\hat{I}^{123}+\hat{I}^{345}+J)\cap R_{1245}$$ 
where $J$ is the {\it right} ideal generated by $\frac{d}{x_{i,3}},~i=1,...,n$ and $R_{1245}$ is the ring of differential operators in $x_{i,1},x_{i,2},x_{4,i},x_{5,i},~i=1,...,n$
Then $\hat{I}^{1245}$ should be symmetric with respect to interchanging of the variables\footnote{Notice that $\hat{I}^{1245}$ is a left ideal in the ring $R_{1245}$.} $x_{i,2}$ and $x_{i,4},~i=1,...,n.$ This property means associativity of our kernel.

{$\bf 3^{\prime}$.} There exists a non-zero solution $K$ of the system (\ref{quang}) is a multiplication kernel of birational type.

Notice that the informal discussion about the properties {\bf P1} and {\bf P2} of the system (\ref{quan}) is also applicable here to the properties {\bf 3} and {$\bf 3^{\prime}$}.

{\bf Remark 3.2.1.}  Weinstein  category (in $C^{\infty}$ setting) discussed in Section 3.1 can be thought of as a semi-classical limit of the symmetric monoidal category of 
complex Hilbert spaces, infinite-dimensional in general. Informally, for a real symplectic manifold $(M,\omega)$ and a small positive Planck constant $\hbar \ll 1$ 
one constructs the corresponding Hilbert space ${\mathcal H}_{\hbar}(M,\omega)$, the quantization of $M$. The dimension of this space $\dim {\mathcal H}_{\hbar}(M,\omega)\approx \frac{1}{(2\pi\hbar)^n}\int_M\frac{\omega^n}{n!}$ where 
$n=\frac{1}{2}\dim M$. Notice that $\dim {\mathcal H}_{\hbar}(M,\omega)$ is infinite if $M$ has infinite volume, for example if $M=T^*X$ is a cotangent bundle. The Hilbert space ${\mathcal H}_{\hbar}(M,-\omega)$ is the dual 
to ${\mathcal H}_{\hbar}(M,\omega)$. The product of symplectic manifolds corresponds to the tensor product of Hilbert spaces. A Lagrangian submanifold $L\subset (M,\omega)$ 
corresponds approximately to a vector $\psi_L\in {\mathcal H}_{\hbar}(M,\omega)$ defined up to multiplication by a phase.  This picture relates collections of commuting operators provided by quantum integrable systems 
 with multiplication kernels (see the beginning of Introduction).

There is another type of quantization of real symplectic manifolds. Namely, with $(M,\omega)$ one can associate an $A_{\infty}$-category depending on small parameter 
$e^{-\frac{1}{\hbar}}$, the Fukaya category ${\cal F}(M,\omega)$. The passing to the opposite manifold $(M, -\omega)$ corresponds to the passing to the opposite 
category, and the product of manifolds corresponds to the tensor product of $A_{\infty}$-categories. Hence, a semi-classical multiplication kernel (i.e. the structure 
of an integrable system with a Lagrangian section) corresponds to a commutative monoids in symmetric monoidal $(\infty,1)$-category of small $A_{\infty}$-categories. In other words, we get a symmetric 
monoidal $A_{\infty}$-category. The basic example of such a category is $\text{Perf}(Y)$ where $Y$ is an algebraic variety endowed with the tensor product over ${\cal O}_Y$. 
In this way one gets a homological mirror symmetry equivalence ${\cal F}(M,\omega) \sim \text{Perf}(Y)$.

\subsection{Example corresponding to Hitchin systems for rank 2 bundles on the projective line with 4 regular singular points}

In this section we describe the semi-classical and quantum multiplication kernels corresponding to Example 2.3.1 in Section 2.3.

Fix an affine coordinate $x$ on $\P^1$, and choose a point on $\P^1$ with coordinate $t\ne 0,1,\infty$. Define $M_0=T^*(\P^1\setminus\{0,1,t,\infty\}$, it is a subset 
of $\A^2=T^*\A^1$ with coordinates $x,p$.

Let $B=\A^1$ and define a fibration $f:~M_0 \to B$ by
$$f(x,p)=x(x-1)(x-t)p^2-s^2x-\frac{k_1^2t}{x}+\frac{k_2^2(t-1)}{x-1}-\frac{k_3^2t(t-1)}{x-t}.$$
Any fiber of $f$ is an elliptic curve (maybe degenerate) with 4 punctures. We define $M\supset M_0$ as a partial compactification obtained by adding 4
missing  points 
on $f^{-1}(b)$ for each $b\in B=\A^1$, i.e. adding 4 copies of $\A^1$. 

Define a Lagrangian submanifold $P\subset (M,-\omega)\times (M,-\omega)\times (M,\omega)$ by
\begin{equation} \label{exc}
\begin{split}
f(x_1,p_1)=f(x_2,p_2)=f(x_3,p_3),~~~~~~~~~~~~~~~~~~~~~~~~~~~~~~~\\
\frac{x_1(x_1-1)(x_1-t)}{(x_1-x_2)(x_1-x_3)}p_1+\frac{x_2(x_2-1)(x_2-t)}{(x_2-x_1)(x_2-x_3)}p_2-\frac{x_3(x_3-1)(x_3-t)}{(x_3-x_1)(x_3-x_2)}p_3-s=0
\end{split}
\end{equation}
where $t,k_1,k_2,k_3,s$ are arbitrary parameters.

Define Lagrangian submanifold $N\subset M$ as the glued copy of $B=\A^1$ corresponding to the puncture $x=\infty$ on the generic fiber of the original map 
$f:M_0=T^*\A^1\to \P^1$.

{\bf Theorem 3.3.1.} The tuple $((M,\omega),N,P)$ constructed above in this Subsection is a semi-classical multiplication kernel.

{\bf Proof.} All properties of a semi-classical multiplication kernel listed in the beginning of this Section can be verified by direct 
computation. In particular, the last equation in the system (\ref{exc}) corresponds to addition on the elliptic curve embedded into $\A^2$ 
with coordinates $x,p$ and defined by the equation
$$x(x-1)(x-t)p^2-s^2x-\frac{k_1^2t}{x}+\frac{k_2^2(t-1)}{x-1}-\frac{k_3^2t(t-1)}{x-t}=b.$$
 $\square$ 
 
 To quantize the above system introdice the differential operators
 \begin{equation}\label{D4}
 D_x=\frac{d}{dx}\cdot x(x-1)(x-t)\cdot\frac{d}{dx}-s(s+2)x-\frac{k_1^2t}{x}+\frac{k_2^2(t-1)}{x-1}-\frac{k_3^2t(t-1)}{x-t},
 \end{equation}
 $$L=\frac{x_1(x_1-1)(x_1-t)}{(x_1-x_2)(x_1-x_3)}\cdot \frac{d}{dx_1}+\frac{x_2(x_2-1)(x_2-t)}{(x_2-x_1)(x_2-x_3)}\cdot \frac{d}{dx_2}+\frac{x_3(x_3-1)(x_3-t)}{(x_3-x_1)(x_3-x_2)}\cdot \frac{d}{dx_3}-s.$$
 
 Consider the following system of equations
 \begin{equation} \label{exq}
\begin{split}
D_{x_1}K(x_1,x_2,x_3)=D_{x_2}K(x_1,x_2,x_3)=D_{x_3}K(x_1,x_2,x_3),\\
LK(x_1,x_2,x_3)=0.~~~~~~~~~~~~~~~~~~~~~~~~~
\end{split}
\end{equation}
 {\bf Theorem 3.3.2.} The system (\ref{exq}) satisfies the properties {\bf 1}, {\bf 2}, {$\bf 3^{\prime}$} listed after the system (\ref{quang}). Moreover, let
 $$K(x_1,x_2,x_3)=(x_1x_2x_3)^sF(u,v)$$
 where
\begin{equation} \label{K}
u=\frac{t^2(x_1-1)(x_2-1)(x_3-1)}{(t-1)^2x_1x_2x_3},~~~v=\frac{(x_1-t)(x_2-t)(x_3-t)}{(t-1)^2x_1x_2x_3},
\end{equation}
and
$$F(u,v)=u^{-k_2}v^{-k_3}\int_{\Gamma}q_1^{2k_2}q_2^{2k_3}\Big(1+q_1+q_2\Big)^{s-k_1-k_2-k_3}\Big(1+\frac{u}{q_1}+\frac{v}{q_2}\Big)^{s+k_1+k_2+k_3}~\frac{dq_1}{q_1}\wedge\frac{dq_2}{q_2}$$
where $\Gamma$ is a cycle. Then $K(x_1,x_2,x_3)$ satisfies the system (\ref{exq}), and if we write 
$$K(x_1,x_2,x_3)=\int_{\Gamma}\tilde{K}(x_1,x_2,x_3,q_1,q_2)~dq_1\wedge dq_2$$
where $\tilde{K}$ is defined by (\ref{K}), then
$$\tilde{K}(x_1,x_2,x_3,q_1,q_2)\tilde{K}(x_3,x_4,x_5,q_3,q_4)~dx_3\wedge dq_1\wedge dq_2\wedge dq_3\wedge dq_4=$$
$$\tilde{K}(x_1,x_4,\tilde{x}_3,\tilde{q}_1,\tilde{q}_2)\tilde{K}(\tilde{x}_3,x_2,x_5,\tilde{q}_3,\tilde{q}_4)~d \tilde{x}_3\wedge d \tilde{q}_1\wedge d \tilde{q}_2\wedge d \tilde{q}_3\wedge d \tilde{q}_4$$
for some birational mapping $(x_3,q_1,q_2,q_3,q_4)\to (\tilde{x}_3,\tilde{q}_1,\tilde{q}_2,\tilde{q}_3,\tilde{q}_4)$.

 {\bf Proof.} The Property {\bf 1} can be verified by direct computation\footnote{In fact, the system (\ref{exq}) was obtained as a quantization of the system (\ref{exc}) with the Properties {\bf 1}, {\bf 2}.} and the Property {\bf 2} is clear because $L$ is symmetric with respect to interchanging of $x_1$ and $x_2$. 
 
 Substituting $K(x_1,x_2,x_3)$ in the form (\ref{K}) in the system (\ref{exq}) one can check that the last equation $LK=0$ holds for an arbitrary function $F(u,v)$ and the first two equations are equivalent to the following system for the function $F(u,v)$
 \begin{equation} \label{G}
\frac{\partial^2G}{\partial u_1\partial v_1}=\frac{\partial^2G}{\partial u_2\partial v_2}=\frac{\partial^2G}{\partial u_3\partial v_3}
\end{equation}
where\footnote{In this formula $k_1,k_2,k_3$ are defined up to sign because the differential operator $D_x$ depends on $k_1^2,k_2^2,k_3^2$ only. Different choice may give different solutions of the system (\ref{exq}) and any solution can be obtained in this way.}
\begin{equation} \label{G1}
G=F\Big(\frac{u_1v_1}{u_3v_3},\frac{u_2v_2}{u_3v_3}\Big)\Big(\frac{v_1}{u_1}\Big)^{k_2}\Big(\frac{v_2}{u_2}\Big)^{k_3}\Big(\frac{v_3}{u_3}\Big)^{k_1}(u_3v_3)^s.
\end{equation}
An Euler type integral representation for a solution of this system can be obtained using the theory of generalized hypergeometric functions. Let us write
$$G=\int_{\Gamma}t_1^{2k_2}t_2^{2k_3}t_3^{2k_1}(u_1t_1+u_2t_2+u_3t_3)^{s-k_1-k_2-k_3}\Big(\frac{v_1}{t_1}+\frac{v_2}{t_2}+\frac{v_3}{t_3}\Big)^{s+k_1+k_2+k_3}\frac{dt_1}{t_1}\wedge\frac{dt_2}{t_2}\wedge\frac{dt_3}{t_3}.$$
It is clear that this expression satisfies (\ref{G}). On the other hand, after change of variables $t_1=\frac{q_1q_3}{u_1},~t_2=\frac{q_2q_3}{u_2},~t_3=\frac{q_3}{u_3}$ and integrating out $q_3$ we obtain the representation (\ref{G1}) where 
$F(u,v)$ are given by (\ref{K}). 

Finally, notice that the family of kernels $\tilde{K}(x_1,x_2,x_3,q_1,q_2)$ is essentially the same as the family of kernels $K(x_1,x_2,y,q_1,q_2)$ in the Example 2.3.1, Section 2.3 (see also Theorem 2.3.1). These kernels are related by a gauge transformation of the form $\tilde{K}(x_1,x_2,x_3,q_1,q_2)\mapsto K(x_1,x_2,x_3,q_1,q_2)\frac{q(x_1)q(x_2)}{q(x_3)}$, a change of variables $(x_1,x_2,x_3,t)\mapsto (\frac{1}{x_1},\frac{1}{x_2},\frac{1}{y},\frac{1}{t})$, and redefinition of other parameters.  $\square$

{\bf Remark 3.3.1.} Differential operators $D_x$ and $L$ in the system (\ref{exq}) can be written in the form 
$$D_x=\frac{d}{dx}\cdot (x^3+g_2x+g_3)\cdot\frac{d}{dx}-s(s+2)x+\frac{k_0+k_1x+k_2x^2}{x^3+g_2x+g_3},$$
$$L=\frac{x_1^3+g_2x_1+g_3}{(x_1-x_2)(x_1-x_3)}\cdot \frac{d}{dx_1}+\frac{x_2^3+g_2x_2+g_3}{(x_2-x_1)(x_2-x_3)}\cdot \frac{d}{dx_2}+\frac{x_3^3+g_2x_3+g_3}{(x_3-x_1)(x_3-x_2)}\cdot \frac{d}{dx_3}-s$$
where $g_2,g_3,k_0,k_1,k_2,s$ are arbitrary constants. This form is convenient to study degenerations of this family. Moreover, if $s=-1$, 
$k_0+k_1x+k_2x^2=(q_1+q_2x)^2$ for some constants $q_1,q_2$, then a solution of the system (\ref{exc}) can be written is a form
$$K(x_1,x_2,x_3)=\exp\Big(q_1f_1(x_1,x_2,x_3)+q_2f_2(x_1,x_2,x_3)\Big)\frac{1}{P(x_1,x_2,x_3)^{-1/2}}$$
where
$$P(x_1,x_2,x_3)=2x_1x_2x_3(x_1+x_2+x_3)-x_1^2x_2^2-x_1^2x_3^2-x_2^2x_3^2+2g_2(x_1x_2+x_1x_3+x_2x_3)+$$
$$4g_3(x_1+x_2+x_3)-g_2^2$$
and $f(x_1,x_2,x_3)=q_1f_1(x_1,x_2,x_3)+q_2f_2(x_1,x_2,x_3)$
satisfies the system of differential equations
$$\frac{\partial f}{\partial x_3}=\frac{1}{(x_3^3+g_2x_3+g_3)P(x_1,x_2,x_3)^{1/2}}\Big(q_1(2x_3^2+x_1x_3+x_2x_3-x_1x_2+g_2)-$$
$$q_2(x_1x_2x_3-x_1x_3^2-x_2x_3^2+g_2x_3+2g_3\Big)$$
and other two equations are obtained from this by a cyclic permutation of $x_1,x_2,x_3$. This family of multiplication kernels coincides, up to a gauge transformation and redefining constants, with a family from the Example 2.2.3, Section 2.2.

{\bf Remark 3.3.2.} Let $D_x$ be differential operator given by (\ref{D4}). Define a function \newline $K_4(x_1,x_2,x_3,x_4)$ by 
$$K_4=F\Big(\frac{x_1x_2x_3x_4}{t^2},\frac{(x_1-1)(x_2-1)(x_3-1)(x_4-1)}{(t-1)^2},\frac{(x_1-t)(x_2-t)(x_3-t)(x_4-t)}{t^2(t-1)^2}\Big)$$
where $F(u,v,w)$ satisfies a system of differential equations
$$\frac{\partial^2G}{\partial u_1\partial v_1}=\frac{\partial^2G}{\partial u_2\partial v_2}=\frac{\partial^2G}{\partial u_3\partial v_3}=\frac{\partial^2G}{\partial u_4\partial v_4}.$$
Here $G$ is given by
$$G=F\Big(\frac{u_2v_2}{u_1v_1},\frac{u_3v_3}{u_1v_1},\frac{u_4v_4}{u_1v_1}\Big)\cdot\Big(\frac{v_2}{u_2}\Big)^{k_1}\Big(\frac{v_3}{u_3}\Big)^{k_2}\Big(\frac{v_4}{u_4}\Big)^{k_3}\Big(\frac{v_1}{u_1}\Big)^{s+1}\cdot\frac{1}{u_1v_1}.$$
Then $K_4$ satisfies the equations
$$D_{x_1}K_4=D_{x_2}K_4=D_{x_3}K_4=D_{x_4}K_4.$$
One can obtain an Euler type integral representation (similar to the one in Theorem 3.3.2) using the theory of generalized hypergeometric functions \cite{GG}. It looks feasible that $K_4$ gives an example of a structure discussed in Section 1.6, problem {\bf 2.} See also Remark 3.4.2 for the generalization of the kernel $K_4$.

\subsection{Example corresponding to Hitchin systems for rank 2 bundles on the projective line with more than 4 regular singular points}

Fix an affine coordinate $x$ on $\P^1$, and choose $n>1$ pairwise distinct points $\P^1$ with coordinates $t_1,...,t_n \ne 0,1,\infty$. Introduce a 
differential operator
$$D_x=\frac{\partial}{\partial x}\cdot x(x-1)(x-t_1)...(x-t_n)\frac{\partial}{\partial x}-s(s+n+1)x^n-\sum_{i=-1}^n \frac{k_{i+2}^2\prod_{j\ne i}(t_i-t_j)}{x-t_i}$$
where $t_{-1}=0,t_0=1$ and $s,k_1,...,k_{n+2}$ are generic parameters. Define furthermore
$$L_{x_1,...,x_{n+1}}=\sum_{i=1}^{n+1}~\frac{1}{\prod_{j\ne i} (x_i-x_j)}~D_{x_i},$$
$$M_{x_1,...,x_{n+1}}=\sum_{i=1}^{n+1}~\frac{x_i(x_i-1)(x_i-t_1)...(x_i-t_n)}{\prod_{j\ne i} (x_i-x_j)}\cdot\frac{\partial}{\partial x_i}-s.$$

{\bf Theorem 3.4.1.} Define the kernel 
$$K_n(x_1,...,x_{2n+1})=\int K(x_1,...,x_{2n+1},q_1,...,q_{n+1})dq_1...dq_{n+1}$$
where the kernel in the r.h.s. is given by (\ref{Kn}). Then the kernel $K_n$ satisfies the following system of holonomic differential equations
\begin{equation}\label{L}
L_{x_{i_1},...,x_{i_{n+1}}}K_n=0,
\end{equation}
\begin{equation}\label{M}
M_{x_{i_1},...,x_{i_{n+1}}}K_n=0
\end{equation}
where $1\leq i_1<...<i_{n+1}\leq 2n+1.$

{\bf Proof.} It is similar to proof of Theorem 3.3.2 from the previous Section. 
 The system (\ref{L}) has the following
 general solution in terms of an arbitrary function $F(w_1,...,w_{n+1})$
 $$K_n=(x_1...~x_{2n+1})^sF(w_0,...,w_n)$$
 where
 $$w_0=\frac{(x_1-1)...(x_{2n+1}-1)~t_1^2...t_n^2}{x_1...~x_{2n+1}~(t_1-1)^2...(t_n-1)^2},~\text{and}~w_i=\frac{(x_1-t_i)...(x_{2n+1}-t_i)\prod_{j\ne i}t_j^2}{x_1...~x_{2n+1}(t_i-1)^2\prod_{j\ne i}(t_i-t_j)^2},~i=1,...,n.$$
 The system (\ref{M}) can be written in terms of $F$ as
 $$\frac{\partial^2G}{\partial u_1\partial v_1}=\frac{\partial^2G}{\partial u_2\partial v_2}=...=\frac{\partial^2G}{\partial u_{n+2}\partial v_{n+2}}$$
 where
 $$G=F\Big(\frac{u_2v_2}{u_1v_1},...,\frac{u_{n+2}v_{n+2}}{u_1v_1}\Big)\cdot\Big(\frac{u_1}{v_1}\Big)^{k_1}...\Big(\frac{u_{n+2}}{v_{n+2}}\Big)^{k_{n+2}}\cdot(u_1v_1)^s.$$
 Our integral representation for the kernel $K$ is the Euler type representation in the theory of generalized hypergeometric functions \cite{GG}. See also the proof of Theorem 3.3.2 and the system (\ref{Dn}).  $\square$

{\bf Remark 3.4.1.} The operator $D_x$ can be also written in the form 
$$D_x=\frac{\partial}{\partial x}\cdot P_{n+2}(x)\frac{\partial}{\partial x}-s(s+n+1)x^n+\frac{Q_{n+1}(x)}{P_{n+2}(x)}$$
where $P_{n+2}(x),~Q_{n+1}(x)$ are arbitrary polynomials in x of degrees $n+2$ and $n+1$ respectively. In this case the operator $L_{x_1,...,x_{n+1}}$ is given by the 
same formula and 
$$M_{x_1,...,x_{n+1}}=\sum_{i=1}^{n+1}~\frac{P_{n+2}(x_i)}{\prod_{j\ne i} (x_i-x_j)}\cdot\frac{\partial}{\partial x_i}-s.$$

{\bf Conjecture 3.4.1} Holonomic cyclic $D$-module given by (\ref{L}), (\ref{M}) defines a generalized product in the sense of Section 4.1. This $D$-module seems to have 
the semi-classical limit described in abstract terms in Section 4.2, and corresponds to Hitchin integrable system for rank 2 bundles with regular singular points at 
$0,1,t_1,...,t_n,\infty$. The fixed point $u_{i_0}$ in notations of Section 4.2 is point $\infty$. 

{\bf Remark 3.4.2.} Define a family of differential operators by
$$D_x=\prod_{i=1}^n(x-t_i)\cdot\Bigg(\frac{\partial^2}{\partial x^2}-\sum_{j=1}^n\frac{b_{j,1}+b_{j,2}-1}{x-t_j}\cdot\frac{\partial}{\partial x}\Bigg)+\sum_{i=1}^n\frac{b_{i,1}b_{i,2}\prod_{j\ne i}(t_i-t_j)}{x-t_i}$$
where we assume $\sum_{i=1}^n(b_{i,1}+b_{i,2})=n-2$. 

The operator $D_x$ has $n$ regular singular points at $x=t_1,...,t_n$ and solutions of the equation $D_xf(x)=0$ near $x=t_i$ have a form $f(x)=(x-t_i)^{b_{i,j}}(1+O(x-t_i))$ for $i=1,...,n$, $j=1,2$. Moreover, any differential operator with these properties and the same symbol as $D_x$ has  form $D_x+\lambda_1+\lambda_2 x+...+
\lambda_{n-3}x^{n-4}$ where $\lambda_1,...,\lambda_{n-3}$ are arbitrary parameters.

Fix an integer $l$ such that $2\leq l\leq n$ and construct a function $K_{n,l}(x_1,...,x_{l+n-4})$ in $n+l-4$ variables as follows
$$K_{n,l}=F\Bigg(q_1\frac{(x_1-t_1)...(x_{l+n-4}-t_1)}{(x_1-t_l)...(x_{l+n-4}-t_l)},...,q_{l-1}\frac{(x_1-t_{l-1})...(x_{l+n-4}-t_{l-1})}{(x_1-t_l)...(x_{l+n-4}-t_l)}\Bigg)\times$$
$$\big((x_1-t_l)...(x_{l+n-4}-t_l)\big)^{-b_{l+1,1}-...-b_{n,1}}\cdot\prod_{l+1\leq k\leq n}\big((x_1-t_k)...(x_{l+n-4}-t_k))^{b_{k,1}}$$
where
$$q_i=\prod_{\substack{1\leq j\leq l-1\\ j\ne i}}\frac{(t_l-t_j)^2}{(t_i-t_j)^2}\cdot \prod_{l+1\leq j\leq n}\frac{t_l-t_j}{t_i-t_j}$$
and $F$ satisfies the system of differential equations\footnote{One can write solutions of this system in terms of Euler type integrals \cite{GG}.} 
$$\frac{\partial^2G}{\partial u_{1,1}\partial u_{1,2}}=\frac{\partial^2G}{\partial u_{2,1}\partial u_{2,2}}=...=\frac{\partial^2G}{\partial u_{l,1}\partial u_{l,2}}.$$
Here 
$$G=F\Bigg(\frac{u_{1,1}u_{1,2}}{u_{l,1}u_{l,2}},...,\frac{u_{l-1,1}u_{l-1,2}}{u_{l,1}u_{l,2}}\Bigg)\cdot\prod_{\substack{1\leq\alpha\leq l\\ 1\leq\beta\leq 2}}u_{\alpha,\beta}^{-b_{\alpha,\beta}}\cdot(u_{l,1}u_{l,2})^{-b_{l+1,1}-...-b_{n,1}}.$$
Notice that if $l=n$, then in the formulas above we have $-b_{l+1,1}-...-b_{n,1}=0$ and $\prod_{l+1\leq j\leq n}=1$. Moreover, in this case the kernel $K_{n,l}$ is a 
function of $2g+2$ variables (here $g=n-3$ is the genus of the spectral curve), and is covariant with respect to the full symmetry group $S_n\ltimes (\Z/2)^n$ 
of the problem. 

Define furthermore
$$L_{x_1,...,x_{n-2}}=\sum_{i=1}^{n-2}~\frac{1}{\prod_{j\ne i} (x_i-x_j)}~D_{x_i}.$$
Then function $K_{n,l}$ satisfies the following system of holonomic differential equations
\begin{equation}\label{Lg}
L_{x_{i_1},...,x_{i_{n-2}}}K_{n,l}=0
\end{equation}
where $1\leq i_1<...<i_{n-2}\leq l+n-4.$

It will be interesting to understand for which $n,l$ the kernel $K_{n,l}$ gives an example of structures discussed in Section 1.6. 

\section{Separation of variables}

\subsection{Families of functions in one variable, and generalized products}

This section can be considered as a continuation of Section 1.4. We will describe the general framework for the Sklyanin method of separation of variables in a broad 
and informal way.

Recall that in Section 1.4 we consider a situation when we have a collection of functions $\psi_{\lambda_1,...,\lambda_n}(x_1,...,x_n)$ in $n$ variables $x_1,...,x_n$
depending on $n$ parameters $\lambda_1,...,\lambda_n$. These functions are the normalized eigenfunctions of a family of commuting operators. Moreover, these functions 
form in a sense a continuous basis of the algebra of functions in $x_1,...,x_n$, thus giving a new product $*$ on functions. 

Now let us consider a different situation: suppose we have a collection $\phi_{\lambda_1,...,\lambda_n}(x)$ of functions in one variable $x$ depending on $n$ parameters 
$\lambda_1,...,\lambda_n$. Then we can construct a collection of functions in $n$ variables $x_1,...,x_n$ (symmetric under permutations) and depending again on $n$
parameters  by 
$$\psi_{\lambda_1,...,\lambda_n}(x_1,...,x_n)=\phi_{\lambda_1,...,\lambda_n}(x_1)...\phi_{\lambda_1,...,\lambda_n}(x_n).$$
Hence one can expect that the functions $\psi_{\lambda_1,...,\lambda_n}$ form a continuous basis of space of $S_n$-invariant functions in $x_1,...,x_n$. 

Conversely, functions $\psi^*_{x_1,...,x_n}(\lambda_1,...,\lambda_n)=\psi_{\lambda_1,...,\lambda_n}(x_1,...,x_n)$ form a continuous basis of the space of functions in 
$\lambda_1,...,\lambda_n$, this is an analogue of the inverse Fourier transform.

{\bf Example 4.1.1.} Let $\phi_{\lambda_1,...,\lambda_n}(x)=e^{\lambda_1 x+\lambda_2 x^2+...+\lambda_n x^n}$. Then $$\psi_{\lambda_1,...,\lambda_n}(x_1,...,x_n)=e^{\lambda_1(x_1+...+x_n)}e^{\lambda_2(x_1^2+...+x_n^2)}...~e^{\lambda_n(x_1^n+...+x_n^n)}.$$ 
Notice that 
$h_1=x_1+...+x_n,...,h_n=x_1^n+...+x_n^n$ are coordinates on $Sym^n\A^1=\A^n$. Hence, we see that the functions  $\psi_{\lambda_1,...,\lambda_n}$ forms the 
continuous basis of Fourier modes. 

For any $x_1,...,x_{n+1}$ consider the function in $\lambda_1,...,\lambda_n$ given by 
$$(\lambda_1,...,\lambda_n)\mapsto \phi_{\lambda_1,...,\lambda_n}(x_1)...\phi_{\lambda_1,...,\lambda_n}(x_{n+1}).$$ Let us expand this function in the continuous
basis $\psi^*_{x_1,...,x_n}(\lambda_1,...,\lambda_n)$. Then we get the generalization of multiplication formulas in Section 1.3.
$$\phi_{\lambda_1,...,\lambda_n}(x_1)...\phi_{\lambda_1,...,\lambda_n}(x_{n+1})=\int K_{n+1,n}(x_1,...,x_{n+1},y_1,...,y_n)\phi_{\lambda_1,...,\lambda_n}(y_1)...\phi_{\lambda_1,...,\lambda_n}(y_n)dy_1...dy_n$$
for some kernel $K_{n+1,n}$ depending on $2n+1$ variables. This kernel satisfies a generalized associativity condition (see also Section 1.5):
$$\int K_{n+1,n}(x_1,...,x_{n+1},z_1,...,z_n)K_{n+1,n}(z_1,...,z_n,x_{n+2},y_1,...,y_n)dz_1...dz_n$$
is symmetric with respect of permutations of $x_1,...,x_{n+2}$. 
Formally, this property implies that the kernel in $3n$ variables given by 
$$K_{2n,n}(x_1,...,x_n,x^{\prime}_1,...,x^{\prime}_n,y_1,...,y_n)=\int K_{n+1,n}(x_1,...,x_n,x^{\prime}_1,z^{(1)}_1,...,z^{(1)}_n)\cdot $$
$$K_{n+1,n}(z_1^{(1)},...,z_n^{(1)},x^{\prime}_2,z_1^{(2)},...,z^{(2)}_n)...K_{n+1,n}(z_1^{(n-1)},...,z_n^{(n-1)},x^{\prime}_n,y_1,...,y_n)\prod_{\substack{1\leq i\leq n\\ 1\leq j\leq n}}dz_j^{(i)}$$
gives a commutative associative product on the space of $S_n$-invariant functions in $n$ variables. In particular, $K_{2n,n}$ is symmetric with respect of permutations 
of $x^{\prime}_1,...,x^{\prime}_n$. 

The above informal considerations can be done formally in the linear algebra framework.

{\bf Definition 4.1.1.} For $n\geq 1$ the generalized commutative associative product on a vector space $V$ is a linear map 
$$\mu_{n+1,n}:V^{\otimes (n+1)}\to V^{\otimes n}$$
which is invariant under $S_{n+1}\times S_n$-action, i.e. it factorises as 
$$V^{\otimes (n+1)}\twoheadrightarrow (V^{\otimes (n+1)})_{S_{n+1}} \to (V^{\otimes n})^{S_n}\hookrightarrow V^{\otimes n},$$
and the map $\mu_{n+2,n}:V^{\otimes (n+2)}\to V^{\otimes n}$ given by $\mu_{n+1,n}\circ (\mu_{n+1,n}\otimes id_V)$ is $S_{n+2}\times S_n$-invariant.
For $n=2$ this map can be represented by the picture
 $$
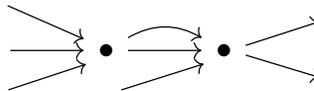
\begin{tikzcd}[row sep=tiny]
  \arrow[dr, end anchor={[yshift=.1ex]}]&  & & \, \\
\,\arrow{r}  &\bullet \arrow[r]\arrow[r, bend left]& \bullet \arrow[r, start anchor={[yshift=0.3ex]}, end anchor={[yshift=1.8ex]north west} ]\arrow[r, end anchor={[yshift=-.5ex]south west} ] &\,\\
\, \arrow[ur, shift right=0.5ex] & \arrow[ur, shift right=0.5ex]& & \,
 \end{tikzcd}$$
 
 This definition of a generalized commutative associative product makes sense in arbitrary symmetric monoidal category. 
 
 {\bf Proposition 4.1.1.} Let us define $\mu_{2n,n}:V^{\otimes 2n}\to V^{\otimes n}$ as the composition of $\mu_{n+1,n}\otimes id_V^{\otimes k}$ for $k=n-1,...,1,0$. For example, for $n=3$ we have 
 $$\mu_{6,3}=\mu_{4,3}\circ(\mu_{4,3}\otimes id_V)\circ(\mu_{4,3}\otimes id_V^{\otimes 2}).$$ 
 This map can be also represented by the picture
  $$
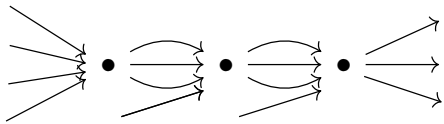
\begin{tikzcd}[row sep=tiny]
  \arrow[dr, start anchor={[yshift=1ex]}]&  &  & & \, \\
\,\arrow[r, start anchor={[xshift=-.2ex,yshift=1.4ex]},
end anchor={[yshift=.1ex]} ] & \bullet \arrow[r]\arrow[r, bend left]\arrow[r, bend right]&\bullet \arrow[r]\arrow[r, bend left]\arrow[r, bend right]&\bullet \arrow{ur}\arrow{r}\arrow[dr, shift right=0.35ex] &\,\\
 \arrow[ur,  start anchor={[yshift=1ex]}] & \arrow[ur, start anchor={[yshift=-1.4ex]}, end anchor={[yshift=-1.4ex]}]  \arrow[ur, start anchor={[yshift=-1.4ex]}, end anchor={[yshift=-1.4ex]}]  &\arrow[ur, start anchor={[yshift=-1.4ex]}, end anchor={[yshift=-1.4ex]}]   & &\,\\
 \arrow[uur, start anchor={[xshift=.3ex,yshift=1ex]}]
 \end{tikzcd}$$
Then $\mu_{2n,n}$ defines a structure of a commutative associative (possibly non-unital) algebra on $Sym^nV$.

\subsection{Semi-classical generalized products}

In this section we describe a class of Poisson compactifications (introduced in Section 8.3 in  \cite{KS}) of cotangent bundles to curves, and a geometric construction 
of semi-classical generalized products. 

Let $C$ be a smooth, not necessarily compact, algebraic curve over $\C$. Write $C=\overline{C}\setminus S$ where $S=\{s_1,...,s_m\}$ is a 
finite subset of the compact curve $\overline{C}$. Denote by ${\mathcal P}_0$ the compact surface which is the total space of $\P^1$-bundle over $\overline{C}$, given by
$\P(\mathbb{O}_{\overline{C}}\oplus T^*_{\overline{C},\log S})$. Surface ${\mathcal P}_0$ carries a natural Poisson structure $\gamma_0$ with the symplectic leaf 
isomorphic to $T^*C$. Poisson structure $\gamma_0$ vanishes on the complement ${\mathcal P}_0\setminus T^*C$, which is the union of smooth divisors 
$\overline{C}_{\infty}\cong \overline{C}$ and $T^*_{s_i}\overline{C}\cong \A^1,~i=1,...,m$ intersecting transversally.

Symplectic form $\omega=\gamma^{-1}_0$ has poles of order 2 along the horizontal divisor $\overline{C}_{\infty}$, and of order one along vertical divisors 
$T^*_{s_i}\overline{C}$. Starting with Poisson surface ${\mathcal P}_0$, let us construct a sequence of Poisson surfaces $({\mathcal P}_i,\gamma_i),~i=0,1,...,N$ 
for some $N\geq 0$ recursively by applying  a sequence of blowups of the following type. Let $p=p_i$ be a point\footnote{Our sequence ${\mathcal P}_i$ depends on choices of these points.} in ${\mathcal P}_i$ such that 
Poisson tensor $\gamma_i$ vanishes at $p$, and there are local coordinates $x,y$ near $p$ such that $\gamma_i=x^ay^b(1+O(x,y))\partial_x\partial_y$, point $p$ has 
coordinates $x=y=0$, and $a+b\geq 2$ for $a,b\in\Z_{\geq0}$. Then we make a blowup at $p$ and obtain a new Poisson surface ${\mathcal P}_{i+1}=Bl_p {\mathcal P}_i$ with 
Poisson tensor $\gamma_{i+1}$. The Poisson tensor $\gamma_{i+1}$ vanishes with order $a+b-1$ at the exceptional divisor. 

Let ${\mathcal P}={\mathcal P}_N$ be the final term of our sequence. The divisor $D$ of zeros of the Poisson structure $\gamma=\gamma_N$ on ${\mathcal P}$ has simple 
normal crossing. The open dense symplectic leaf ${\mathcal P} \setminus D$ is equal to $T^*C$. 

Denote by $\{D_\alpha\}$ the set of irreducible components of  $D$ at which $\gamma$ vanishes with multiplicity one (or, equivalently, the 
meromorphic symplectic form $\omega=\gamma^{-1}$ has pole of order one)\footnote{We will not use other components of $D$ in our considerations.}. It follows by induction from the construction that divisors  $D_{\alpha}$  do not intersect each other. Moreover, each $D_{\alpha}$ contains exactly one double point 
of $D$, and the complement $D^0_{\alpha}$ to this point is isomorphic to $\A^1$. In fact, there is a canonical coordinate on $D^0_{\alpha}$ given by the residue of the 
restriction of the Liouville 1-form on $T^*C$ to a small disc in ${\mathcal P}$ transversally intersecting  $D^0_{\alpha}$.

{\bf Remark 4.2.1.} One can associate with each component $D_{\alpha}$ as above, a point $v\in\overline{C}$ and a Puiseaux series 
$f(x)\in \cup_{n\geq 1}\C(\!(x^{\frac{1}{n}})\!)=\overline{\C(\!(x)\!)}$ where $x$ is a local coordinate at $v$, up to certain identification. First, we identify two series 
differ by $\cup_{n\geq 1}\C[[x^{\frac{1}{n}}]]$. In other words, we can keep only terms with strictly negative exponents. Second, we identify series which differ by the action 
of the Galois group $\widehat{\Z}$ of $\overline{\C(\!(x)\!)}$, with the topological generator $x^{\frac{1}{n}}\mapsto e^{\frac{2\pi i}{n}}x^{\frac{1}{n}}$. 

For a given component $D_{\alpha}$ the corresponding Puiseaux series is defined as follows. Let us choose a germ $\Sigma_{\alpha}$ of smooth curve in ${\mathcal P}$ 
intersecting $D^0_{\alpha}$ transversely at one point, and such that $\Sigma_{\alpha}$ projects to $\overline{C}$ non-trivially, i.e. not to a point. Then the punctured 
disc $\Sigma_{\alpha}\setminus \Sigma_{\alpha}\cap D^0_{\alpha}$ contained in $T^*C$ can be considered as the graph of a meromorphic multivalued 1-form at a point 
$v\in\overline{C}$. This form can be written in local coordinate $x$ at $v$ as $df(x)+\lambda d \log x$ where $\lambda\in\C,~f(x)\in\overline{\C(\!(x)\!)}$. One can show 
that different choices of germ $\Sigma_{\alpha}$ give equivalent series $f(x)$ in the above sense. Incidentally, the equivalence classes of such series correspond to 
all the possible irregular terms for formal meromorphic connections on $\C(\!(x)\!)$. 

Let us choose a collection of points $\sum_i n_iu_i$, $n_i\geq 1$ in $\coprod_{\alpha}D^0_{\alpha}=\coprod_{\alpha}\A^1_{\alpha}$ with multiplicities, and an integer 
$r\geq 1$. With the tuple $({\mathcal P},\sum_i n_i u_i,r)$ we associate a classical integrable system. The base $B$ of this system will be the set of smooth connected 
curves $\Sigma\subset{\mathcal P}$ (spectral curves) such that $\Sigma\not\subset D$ and $\Sigma\cap D=\sum_i n_i u_i$, and the projection $\Sigma
\to \overline{C}$ has degree $r$. 

Let us assume that $B\ne \emptyset$ (this assumption considered as a property of $({\mathcal P},\sum_i n_i u_i,r)$ is related with the additive Deligne-Simpson problem 
\cite{S}). One can show that $B$ is 
an open dense subset of $\A^g$ where $g$ is the genus of any spectral curve $\Sigma$. The tangent space $T_{\Sigma}B$ is canonically identified with
 $\Gamma(\Sigma,\Omega^1_{\Sigma})$. 
 
 Let us fix an integer $d\in\Z$. Define $M_d$ to be the space of pairs $(\Sigma,[L])$ where $[L]\in Pic_d(\Sigma)$ is a class of a line bundle of degree $d$. 
 
 There is a natural symplectic form $\omega_{M_d}$ on $M_d$, and the natural projection $M_d\to B$ is a Liouville integrable system.
 
 The above construction gives an alternative description of a Zariski open dense part of Hitchin integrable systems for group $GL_r$ on $\overline{C}$ with possibly irregular 
 singularities.
 
 Notice that among integrable systems $M_d\to B$ only one, corresponding to $d=0$, has an obvious Lagrangian section, which is the zero section. 
 
 On the other hand, we can define a birational symplectomorphism\footnote{We use the notation $\sim$ for birational equivalence.} $M_g\sim Sym^gT^*C\sim T^* Sym^gC.$ Indeed, we have $Pic_g(\Sigma)\sim Sym^g \Sigma$. Hence, a generic point 
 in $M_g$ is a spectral curve $\Sigma$ and a collection of $g$ points $(t_1,...,t_g)$ of $\Sigma$ up to a permutation. Generically all points $t_i$ will be distinct and 
 do not belong to $\Sigma\cap D$, hence give $g$ points in $T^*C$. Conversely, given $g$ generic points in $T^*C$, there exists a unique spectral curve $\Sigma$ passing 
 through them, and therefore we obtain a point in $M_g$. 
 
 This is a geometric version of method of separation of variables. This is different from the usual Hitchin systems, where the total space is identified birationally with 
 the cotangent bundle to the moduli space $Bun_G(\overline{C})$ of $G$-bundles on $\overline{C}$, or its version associated with marked points and irregular 
 singularities.

Now we can define a semi-classical generalized product. Let us pick a point $u_{i_0}$ among the collection $\{u_i\}$. Notice that $u_{i_0}\in\Sigma$ for any spectral 
curve $\Sigma\in B$. Using $u_{i_0}$ we can identify $M_d$ for all $d\in \Z$ by adding multiples of $u_{i_0}$. Also, we define a Lagrangian subvariety 
$$L_{g+1,g}\subset Sym^{g+1}(T^*C,-\omega)\times Sym^g(T^* C,\omega)$$
by the formula 

$L_{g+1,g}=\{(a_1,...,a_{g+1},b_1,...b_g);$ there exists $\Sigma\in B$ such that 
$$ a_1,...,a_{g+1},b_1,...b_g\in\Sigma,~~~a_1+...+a_{g+1}=b_1+...+b_g+u_{i_0}\in Pic(\Sigma)\}.$$
This variety might be singular, so we treat it only as a first rough approximation.

It is clear that $L_{g+1,g}$ considered as a morphism in Weinstein category, is a generalized product. 

Similarly, one can define Lagrangian correspondences $L_{g+k,g}$ for any $k\geq 1$. In the case $k=g$ we obtain a semi-classical multiplication kernel for the symplectic 
variety $M_0=M_g$. 

{\bf Remark 4.2.2.} The choice of a point $u_{i_0}$ corresponds in a sense to the choice of a normalization of functions $\phi_{\lambda_1,...,\lambda_n}(x)$ in Section 4.1, and therefore, to the choice of a normalization of functions $\psi_{\lambda_1,...,\lambda_n}(x_1,...,x_n)$.

\subsection{Towards quantized generalized product}

In this section we will propose a hypothetical construction of the quantization of the semi-classical generalized product introduced in Section 4.2, understood as a 
holonomic $D$-module on $C^{2g+1}$ together with a cyclic vector. 

Let $({\mathcal P},\sum_i n_iu_i,r)$ be a tuple as in Section 4.2. We now assume that for any two different points $u_i\ne u_j$ belonging to the same component $D^0_{\alpha}=\A^1$ 
the difference $\mu_i-\mu_j$ between their canonical coordinates is not an integer. 
Starting with such a tuple we construct a family of cyclic $D$-modules on $\overline{C}$ depending on $g$ parameters $\lambda_1,...,\lambda_g$. 
Namely, consider meromorphic differential operators $L$ of order $r$ acting from the trivial line bundle ${\cal O}_{\overline{C}}$ to $K^{\otimes r}_{\overline{C}}$, 
with symbol $\Big(\frac{\partial}{\partial x}\Big)^r$ in local coordinates, and such that singularities of solutions correspond to $\sum_i n_i u_i$. To explain the latter 
condition more precisely, recall that a point $u_i\in D$ belongs to a component $D_{\alpha_i}$ of the divisor $D$. In Section 4.2 we explained that divisor $D_{\alpha_i}$
 gives a Puiseaux series $f_i(x)$ at point $v_i=pr_{D\to\overline{C}}(u_i)\in \overline{C}$. Denote by $x=x_i$ the local coordinate at $v_i$, and by $\mu_i\in\C$ the 
position of point $u_i\in D$ on $D^0_{\alpha_i}=\A^1$. Then we say that differential equation $L\phi=0$ has singularity at $v_i$ corresponding to $u_i$ and $n_i\geq 1$ 
if this equation has solutions with asymptotic behaviour 
$$\phi(x)=x^{\mu_i}e^{f_i(x)}\log(x)^k(1+...),~~~k=0,...,n_i-1.$$
One can show that differential operators $L$ satisfying these properties form an affine space of dimension $g$ (a version of variety of opers for the Hitchin system 
for group $GL_r$), and can be written as $L=L_0+\sum_{i=1}^g\lambda_i L_i$ where $\deg L_0=r,~\deg L_i<r,i=1,...,g$. Here $(\lambda_1,...,\lambda_g)$ are parameters 
(coordinates on the space of opers). 

Recall that in Section 4.2 we made a choice $u_{i_0}$ of one point of $\Sigma\cap D$. Let us assume that the corresponding Puiseaux series $f_{i_0}$ is unramified
\footnote{We do not know how to extend our construction to the case of ramified $f_{i_0}$.}, i.e. belong to $\C(\!(x)\!)\subset \cup_{n\geq1}\C(\!(x^{\frac{1}{n}})\!)$. We define the
 normalized solution corresponding to $u_{i_0}$ as the unique formal solution at $v_{i_0}$ of the form 
$$\phi_{\lambda_1,...,\lambda_g}(x)=x^{\mu_{i_0}}e^{f_{i_0}(x)}(1+\sum_{j\geq 1}P_jx^j)$$
where $P_1,P_2,...$ depend on $\lambda_1,...,\lambda_g$. One can see that $P_j$ are polynomials in $\lambda_1,...,\lambda_g$, we set $P_0=1$. 

Without loss of generality, in order to simplify the exposition, we assume that $\mu_{i_0}=0,~f_{i_0}(x)=0$. This can be achieved by the conjugation of $L$. 

Examples suggest that the following is true:

\framebox{The set $\{P_{i_1}...P_{i_g};~i_1\leq...\leq i_g\}$ forms a linear basis of $\C[\lambda_1,...,\lambda_g]$.}

Assuming this property, we can define coefficients 
$$C^{j_1,...,j_g}_{i_1,...,i_{g+1}}=\text{Coeff}_{P_{j_1}...P_{j_g}}(P_{i_1}...P_{i_{g+1}}).$$
These coefficients are structure constants of a generalized commutative associative product on an infinite-dimensional space with a basis $e_0,e_1,...$ where $e_j$ 
corresponds to $P_j$. 

{\bf Remark 4.3.1.} In general, suppose that  $A$ is a commutative associative algebra over a field $k$ of characteristic zero, and $V\subset A$ is a vector subspace such that the composition $Sym^g V\to V^{\otimes g}\to A^{\otimes g}\to A$ is an isomorphism of vector spaces, where the last map is induced by multiplication in $A$. Then $V$ 
carries a structure of commutative associative generalised product $Sym^{g+1} V\to Sym^g V=A$. 

Let us consider generating series 
\begin{equation}  \label{2}
K_{g+1,g}(x_1,...,x_{g+1},y_1,...,y_g)=\sum_{i_1,...j_g\geq 0}C^{j_1,...,j_g}_{i_1,...,i_{g+1}}x_1^{i_1}...x_{g+1}^{i_{g+1}}y_1^{-j_1}...y_g^{-j_g}\prod_{k=1}^g\frac{dy_k}{y_k}.
\end{equation}
We consider $K_{g+1,g}$ as an element of $W/W_+$ where $W=\C[[y_1,...,y_g]][y_1^{-1},...,y_g^{-1}][[x_1,...,x_{g+1}]]$ and $W_+\subset W$ is the subspace consisting of 
series which do not have poles in variable $y_i$ for some $i=1,...,g$. 

Consider algebraic variety $X=\overline{C}^{2g+1}$ with the line bundle ${\mathcal L}=\otimes^{2g+1}_{i=g+2}pr^*_iT^*_{\overline{C}}$ and a point 
$p=\underbrace{(v_{i_0},...,v_{i_0})}_{2g+1\text{~times}}$. Then we have an algebra $\text{Dif}_{rat}$ of differential operators on ${\mathcal L}$ with coefficients 
in rational functions on $X$ (i.e. differential operators at the generic point of $X$). It has a subalgebra $\text{Dif}_{rat,p}$ of differential operators without 
poles at $p$. Vector spaces $W,W_+$ and hence $W/W_+$ are $\text{Dif}_{rat,p}$-modules. 

{\bf Conjecture 4.3.1.} The cyclic $\text{Dif}_{rat,p}$-module ${\sf K}_{g+1,g}$  generated by $K_{g+1,g}$ is holonomic. Moreover, 
$\text{Dif}_{rat}\otimes_{\text{Dif}_{rat,p}}{\sf K}_{g+1,g}$ is the kernel of the generalized product localized at the generic point of $X=\overline{C}^{2g+1}$. 

Examples suggest that this kernel can be written in an explicit form, see equations (\ref{kerint}), (\ref{elfun}) in Introduction where the integration is understood 
as the direct image of holonomic $D$-modules. 

{\bf Remark 4.3.2.} The example of a family of functions $\phi_{\lambda_1,...,\lambda_g}(x)=e^{\lambda_1 x+...+\lambda_g x^g}$ does not fit to the above scheme. Here 
$\phi_{\lambda_1,...,\lambda_g}(x)$ is the unique solution of the equation $(L_0+\lambda_1 L_1+...+\lambda_g L_g)\phi(x)=0$ in $\C[[x]]$ with constant term 1. 
Here $L_0=\frac{\partial}{\partial x},~L_i=-ix^{i-1},~i=1,...,g$. It will be interesting to find a general framework which includes this example. 

Conjecture 4.3.1 is formulated in terms of the quotient $W/W_+$ which is a complicated object not suitable for explicit computations. It is more convenient to see $K_{g+1,g}$ 
in (\ref{2}) literally as a function of $x_i,y_i$. For example, we can assume that $\overline{C}=\P^1$ and $x$ is a global coordinate. 

The examples studied in the next section indicate that the cyclic $D$-module for the lifted kernel is still holonomic, and it maps epimorphically to the expected 
generalized multiplication kernel.

{\bf Remark 4.3.3.} The quantum generalized multiplication kernel should produce isomorphisms of the following type. For all $\vec{\lambda}=(\lambda_1,...,\lambda_g)$ we 
expect that 
\begin{equation}\label{3}
(C^{2g+1}\to C^{g+1})_*~({\sf K}_{g+1,g}\otimes (\underbrace{{\cal E}_{\vec{\lambda}}\boxtimes ...\boxtimes {\cal E}_{\vec{\lambda}}}_{g~\text{times}}))\simeq \underbrace{{\cal E}_{\vec{\lambda}}\boxtimes ...\boxtimes {\cal E}_{\vec{\lambda}}}_{g+1~\text{times}}
\end{equation}
where ${\cal E}_{\vec{\lambda}}$ is a cyclic holonomic $D$-module (oper) parameterized by $\vec{\lambda}$. It is known in the theory of integrable systems that the 
locus of opers is a Lagrangian subvariety in the symplectic manifold parameterizing all (non-cyclic) holonomic $D$-modules on $C$ with given singularities (de Rham 
moduli space). We claim that the equation (\ref{3}) holds also for more general modules. The rough reason is that (\ref{3}) makes sense in Betti realization, in which 
the locus of opers is Zariski dense.

\subsection{A generalization to other Poisson surfaces}

In Section 4.2 we made a sequence of blowups of the Poisson surface ${\mathcal P}_0=\P(\mathbb{O}_{\overline{C}}\oplus T^*_{\overline{C},\log S})$. Assume that 
$S\ne\emptyset$, hence $C$ is affine. The quantization of the open dense symplectic leaf of ${\mathcal P}_0$ is the algebra of differential operators on $C$. 
Let us refer to this class of Poisson surfaces as to ``rational''. There are two other classes of compact Poisson surfaces, which we will call ``trigonometric'' and 
``elliptic'', with an open dense symplectic leaf $M$. In all three cases the open dense symplectic leaf $M$ is an affine variety. 

In the trigonometric case ${\mathcal P}$ is any toric compactification of its symplectic leaf $M=\C^*\times\C^*$ with $\omega=\frac{dx}{x}\wedge \frac{dy}{y}$. The 
quantization of $M$ is quantum torus $A_q=\C\langle X^{\pm 1},Y^{\pm 1}\rangle /XY=qYX$. 

In the elliptic case ${\mathcal P}=\C P^2$ with the symplectic leaf $M={\mathcal P}\setminus \{\text{cubic curve}\}$. The quantization is an inhomogeneous version of 
Sklyanin algebra with three generators. 

In the elliptic case, there is also a version with non-affine $M$ similar to $T^*C$ for compact curve $C$. Namely, let ${\mathcal P}=\P^1\times E$ where $E$ is an 
elliptic curve and $M=\C^*\times E$. In this case instead of algebras one should consider abelian or triangulated categories (analogs of categories of modules). 

In our considerations in Section 4.2 the principal role was played by divisors on the blowup at which $\omega$ has pole of order one. 

In the trigonometric case such divisors do not have continuous parameters and correspond to pairs of coprime integers $(a,b)$, i.e. these divisors are toric divisors. 
Each $D^0_{\alpha}$ is isomorphic to $\C^*$. 

In the elliptic case there is only one such a divisor, the initial cubic curve $E$ (and there are no non-trivial blowups). 

Similarly to Section 4.2, we choose a sequence of blowups (trivial in the elliptic case) and a collection of points with multiplicities on 
$\coprod_{\alpha} D^0_{\alpha}$. 

In this way we obtain a classical integrable system. Choosing one point $u_{i_0}$ we get a semi-classical kernel $K_{g+1,g}$. The analog of opers will be cyclic 
holonomic $A$-modules where $A$ is the quantization of $M$. The quantum kernel should be a cyclic holonomic $A^{\otimes (g+1)}\otimes (A^{op})^{\otimes g}$-module. 

Notice that the rational case $M=T^*C$ is related to the geometric Langlands correspondence for groups $GL_r$. In the trigonometric (resp. elliptic) cases the
 multiplication kernels $K_{g+1,g}$ are also expected to be a kind of motivic in trigonometric (resp. elliptic) sense. For example, the $q$-analog of motivic 
 holonomic modules are discussed in Section 6.1 of \cite{K1}. Roughly, these modules are build from the basic $A_q$-modules with cyclic vector $A_q/A_q\cdot (X+Y-1)$ by external 
 tensor products, actions of $Sp(2n,\Z)$ and pushforwards.

\section{Multiplication kernels associated with differential operators}

In this section we always work in the global coordinate on $\A^1$. 

\subsection{General setup}

Let $P_0(\lambda),P_1(\lambda),...\in\C[\lambda]$ be a basis of the vector space $\C[\lambda]$ such that $\deg P_i=i$ and $P_0=1$. We have 
\begin{equation}\label{strc}
P_i(\lambda)P_j(\lambda)=\sum_{k=0}^{\infty}C_{i,j}^k P_k(\lambda)
\end{equation}
where $C_{i,j}^k$ are structure constants of polynomial multiplication in the basis $P_i(\lambda)$. Here we assume that $C_{i,j}^k$ are independent of $\lambda$.

Construct generating functions:
\begin{equation}\label{genf}
f_{\lambda}(x)=\sum_{i=0}^{\infty}P_i(\lambda)x^i,
\end{equation}
\begin{equation}\label{genK}
K(x_1,x_2,y)=\sum_{i,j,k\geq 0}C_{i,j}^k\frac{x_1^ix_2^j}{y^{k+1}}.
\end{equation}
Notice that $f_{\lambda}(x)=1+O(x)$ and $K(x,0,y)=\frac{1}{y-x}$.

We have by construction
\begin{equation}\label{fK}
f_{\lambda}(x_1)f_{\lambda}(x_2)=\frac{1}{2\pi i}\oint K(x_1,x_2,y)f_{\lambda}(y)dy
\end{equation}
where integral is taken by a small circle around zero. This means that the associativity condition holds
\begin{equation}\label{asscont}
\oint K(x_1,x_2,y)K(y,x_3,z)dy=\oint K(x_1,x_3,y)K(y,x_2,z)dy.
\end{equation}

Let $D_x$ be a differential operator in $x$. Assume that our generating function is a solution of the differential equation 
$$D_xf_{\lambda}(x)=\lambda f_{\lambda}(x).$$
In this case the kernel $K(x_1,x_2,y)$ satisfies the equations
$$D_{x_1}K(x_1,x_2,y)=D_{x_2}K(x_1,x_2,y)$$
and 
$$D_{x_1}K(x_1,x_2,y)-D_{y}K(x_1,x_2,y)\in \C[[y]].$$

{\bf Example 5.1.1.} Let $D_x=\frac{d}{dx}$ and $f_{\lambda}(x)=e^{\lambda x}$. In this case we have 
$$K(x_1,x_2,y)=\frac{1}{y-x_1-x_2}.$$

The construction above can be generalized to the case of polynomials $P_i(\lambda_1,...\lambda_g),i=0,1,...$ such that $P_0(\lambda_1,...\lambda_g)=1$ and
$\{P_{i_1}...P_{i_g},~0\leq i_1\leq i_2\leq...\leq i_g\}$ is a basis in the vector space $\C[\lambda_1,...,\lambda_g]$.  We have 
\begin{equation}\label{strcm}
P_{i_1}...P_{i_{g+1}}=\sum_{j_1,...,j_g\geq 0}C_{i_1,...,i_{g+1}}^{j_1,...,j_g}P_{j_1}...P_{j_g}
\end{equation}
where we assume that $C_{i_1,...,i_{g+1}}^{j_1,...,j_g}$  are independent of $\lambda_1,...,\lambda_g$ and symmetric with respect to indexes $j_1,...,j_g$.

Construct generating functions:
\begin{equation}\label{genfm}
f_{\vec{\lambda}}(x)=\sum_{i=0}^{\infty}P_i(\lambda_1,...,\lambda_g)x^i,~~~~\text{where}~~~\vec{\lambda}=(\lambda_1,...,\lambda_g),
\end{equation}
\begin{equation}\label{genKm}
K(x_1,...,x_{g+1},y_1,...,y_g)=\sum_{i_1,...,j_g\geq 0}C_{i_1,...,i_{g+1}}^{j_1,...,j_g}x_1^{i_1}...x_{g+1}^{i_{g+1}}y_1^{-j_1-1}...y_g^{-j_g-1}.
\end{equation}
Notice that $f_{\vec{\lambda}}(x)=1+O(x)$ and $$K(x_1,...,x_n,0,y_1,...,y_n)=\frac{1}{n!}\sum_{\sigma\in S_n}\frac{1}{(y_{\sigma_1}-x_1)...(y_{\sigma_n}-x_n)}.$$

We have by construction 
\begin{equation}\label{relfKm}
f_{\vec{\lambda}}(x_1)...f_{\vec{\lambda}}(x_{n+1})=\frac{1}{(2\pi i)^n}\oint...\oint K(x_1,...,x_{n+1},y_1,...,y_n) f_{\vec{\lambda}}(y_1)...f_{\vec{\lambda}}(y_n)dy_1...dy_n
\end{equation}
where we integrate over small circles around zero with respect to each $y_i$. We also have an associativity condition: the expression
\begin{equation}\label{asscontm}
\oint...\oint K(x_1,...,x_{n+1},y_1,...,y_n) K(x_{n+2},y_1,...,y_n,z_1,...,z_n)dy_1...dy_n
\end{equation}
is symmetric with respect to $x_1,...,x_{n+2}$.

Assume that\footnote{In the notations of Section 4.3 we have $L_0=D_x,~L_i=-x^{i-1},~i=1,...,g$.} 
$$D_xf_{\vec{\lambda}}(x)=(\lambda_1+\lambda_2x+...+\lambda_gx^{g-1}) f_{\vec{\lambda}}(x)$$
for a differential operator $D_x$. In this case the kernel $K(x_1,...,x_{n+1},y_1,...,y_n)$ satisfies the equation
$$\sum_{i=1}^{n+1} \frac{1}{(x_1-x_i)...\hat{i}...(x_{n+1}-x_i)}D_{x_i}K(x_1,...,x_{n+1},y_1,...,y_n)=0.$$

 \subsection{A kernel associated with first order differential operators}
 
 Here we return to the basic example introduced in Section 4.1.
 
Define polynomials $P_i(\lambda_1,...\lambda_g)$ by 
$$e^{\lambda_1x+\lambda_2x^2+...+\lambda_gx^g}=\sum_{i\geq 0}P_i(\lambda_1,...\lambda_g)x^i.$$
Notice that if 
$$D_x=\frac{d}{dx}-(\lambda_1+2\lambda_2x+...+g\lambda_gx^{g-1}),$$
then $D_xe^{\lambda_1x+\lambda_2x^2+...+\lambda_gx^g}=0$.

Define structure constants $C_{i_1,...,i_{g+1}}^{j_1,...,j_g}\in\Q$ by
(\ref{strcm}).

Define kernel $K(x_1,...,x_{g+1},y_1,...,y_g)$ as the generating function   for these structure constants by (\ref{genKm}).

{\bf Theorem 5.2.1.} 
$$K(x_1,...,x_{g+1},y_1,...,y_g)=\frac{1}{g!}\sum_{\sigma\in S_g}\frac{1}{(y_{\sigma_1}-q_1)...(y_{\sigma_g}-q_g)}$$
where $q_1,...,q_g$ are roots of a polynomials
$$Q(t)=t^g-(x_1+...+x_{g+1})t^{g-1}+(x_1x_2+...+x_gx_{g+1})t^{g-2}+...\pm (x_1...x_g+...+x_2...x_{g+1}).$$
Coefficients of this polynomial in $t$ are elementary symmetric polynomials in $x_1,...,x_{g+1}$.

{\bf Proof.} To simplify notations let $g=2$, the general case is similar. We have
$$e^{\lambda_1(x_1+x_2+x_3)+\lambda_2(x_1^2+x_2^2+x_3^2)}=\sum_{i_1,i_2,i_3\geq 0}P_{i_1}P_{i_2}P_{i_3}x_1^{i_1}x_2^{i_2}x_3^{i_3}=
\sum_{i_1,...,j_2\geq 0}C_{i_1i_2i_3}^{j_1j_2}x_1^{i_1}x_2^{i_2}x_3^{i_3}P_{j_1}P_{j_2}=
$$
$$=\frac{1}{(2\pi i)^2}\sum_{i_1,...,j_2\geq 0}C_{i_1i_2i_3}^{j_1j_2}x_1^{i_1}x_2^{i_2}x_3^{i_3}\oint \oint \frac{e^{\lambda_1(y_1+y_2)+\lambda_2(y_1^2+y_2^2)}dy_1dy_2}{y_1^{j_1+1}y_2^{j_2+1}}=$$
$$=
\frac{1}{(2\pi i)^2}\oint\oint K(x_1,x_2,x_3,y_1,y_2)e^{\lambda_1(y_1+y_2)+\lambda_2(y_1^2+y_2^2)}dy_1dy_2$$
where integrals are taken over small circles around $0$. Expanding in power series in $\lambda_1,~\lambda_2$ and equating coefficients at $\lambda_1^{n_1}\lambda_2^{n_2}$ we get
$$(x_1+x_2+x_3)^{n_1}(x_1^2+x_2^2+x_3^2)^{n_2}=\frac{1}{(2\pi i)^2}\oint\oint K(x_1,x_2,x_3,y_1,y_2)(y_1+y_2)^{n_1}(y_1^2+y_2^2)^{n_2}dy_1dy_2.$$
It follows that 
$$\phi(x_1+x_2+x_3,x_1^2+x_2^2+x_3^2)=\frac{1}{(2\pi i)^2}\oint\oint K(x_1,x_2,x_3,y_1,y_2)\phi(y_1+y_2,y_1^2+y_2^2)dy_1dy_2$$
where $\phi(u,v)$ is an arbitrary function analytic near $u=v=0$. We can write this condition as 
$$\phi(x_1+x_2+x_3,x_1x_2+x_2x_3+x_1x_3)=\frac{1}{(2\pi i)^2}\oint\oint K(x_1,x_2,x_3,y_1,y_2)\phi(y_1+y_2,y_1y_2)dy_1dy_2$$
by changing variables in $\phi$ (because elementary symmetric functions can be written in terms of sums of powers). 

Let $\phi=\phi_{m_1,m_2}$ where $\phi_{m_1,m_2}(y_1+y_2,y_1y_2)=y_1^{m_1}y_2^{m_2}+y_2^{m_1}y_1^{m_2}$, it is clear that 
 $\phi_{m_1,m_2}(u,v)=s_1^{m_1}s_2^{m_2}+s_2^{m_1}s_1^{m_2}$ where $s_1,s_2$ are roots of polynomial $t^2-ut+v$. With this choice of $\phi$ we obtain
 $$\frac{1}{(2\pi i)^2}\oint\oint K(x_1,x_2,x_3,y_1,y_2)(y_1^{m_1}y_2^{m_2}+y_2^{m_1}y_1^{m_2})dy_1dy_2=q_1^{m_1}q_2^{m_2}+q_2^{m_1}q_1^{m_2}$$
where $q_1,~q_2$ are roots of polynomial $t^2-(x_1+x_2+x_3)t+x_1x_2+x_2x_3+x_1x_3$. It follows 
$$K(x_1,x_2,x_3,y_1,y_2)=\frac{1}{2}\sum_{m_1,m_2\geq 0}(q_1^{m_1}q_2^{m_2}y_1^{-m_1-1}y_2^{-m_2-1}+q_2^{m_1}q_1^{m_2}y_1^{-m_1-1}y_2^{-m_2-1})$$
and summing up geometric series we obtain the statement of the Theorem. $\Box$

\subsection{The case of second order differential operators with 4 regular singular points}

Let 
$$D_x=x(x-1)(x-t)\frac{d^2}{dx^2}+x(x-1)(x-t)\Big(\frac{s_1}{x}+\frac{s_2}{x-1}+\frac{s_3}{x-t}\Big)\frac{d}{dx}+r_1r_2x+\lambda$$
where $t,s_1,s_2,s_3,r_1,r_2,\lambda$ are parameters such that $t\ne 0,1$ and 
$$r_1+r_2=s_1+s_2+s_3-1.$$
This is the most general second order differential operators with regular singularities at $x=0,1,t,\infty$ and with analytic solutions near $x=0,1,t$.

There exists a unique solution $f_{\lambda}(x)$ of the differential equation 
$$D_xf_{\lambda}(x)=0$$
such that $f_{\lambda}(x)$ is analytic near $x=0$ and $f_{\lambda}(0)=1$. We have
$$f_{\lambda}(x)=\sum_{i=0}^{\infty}P_i(\lambda)x^i$$
where $P_i(\lambda)$ are polynomials in $\lambda$ of degree $i$ and $P_0(\lambda)=1$.

Since $P_i(\lambda)$, $i=0,1,...$ is a basis of the vector space $\C[\lambda]$ we can define structure constants of polynomial multiplication in this basis by (\ref{strc}).

Define kernel $K(x_1,x_2,y)$ as a generating function of these structure constants by (\ref{genK}).

Recall that  the Gauss hypergeometric function is given by
$$F(a,b,c,u)=\sum_{n=0}^{\infty}\frac{a(a+1)...(a+n-1)\cdot b(b+1)...(b+n-1)}{c(c+1)...(c+n-1)}\cdot \frac{u^n}{n!}=$$
$$=\frac{\Gamma(c)}{\Gamma(b)\Gamma(c-b)}\int_0^1p^{b-1}(1-p)^{c-b-1}(1-pu)^{-a}~dp.$$

{\bf Theorem 5.3.1.} The kernel $K(x_1,x_2,y)$ is given by
$$K(x_1,x_2,y)=\frac{1}{y-1}\sum_{i=0}^{\infty}\frac{(r_1+r_2-s_1-s_3+1)...(r_1+r_2-s_1-s_3+i)}{(r_1+r_2-s_1+i)...(r_1+r_2-s_1+2i-1)}\times$$
$$F\Big(i+1,r_1+r_2-s_1-s_3+i+1,r_1+r_2-s_1+2i+1,\frac{t-1}{y-1}\Big)\times$$
$$F\Big(-i,r_1+r_2-s_1+i,r_1+r_2-s_1-s_3+1,\frac{t(x_1-1)(x_2-1)}{(t-1)(x_1x_2-t)}\Big)\times$$
$$F\Big(r_1+i,r_2+i,s_1,\frac{x_1x_2}{t}\Big)~\Big(\frac{x_1x_2}{t}-1\Big)^i~\frac{(t-1)^i}{(y-1)^i}.$$
The proof is based on the following

{\bf Lemma 5.3.1.} The function $K(x_1,x_2,y)$ is the unique function characterized by the following properties:

{\bf 1.} $K(x_1,x_2,y)=K(x_2,x_1,y)$.

{\bf 2.} $K(x_1,x_2,y)$ has Laurent series expansion by non negative powers of $x_1,x_2$ and negative powers of $y$.

{\bf 3.} $D_{x_1}K(x_1,x_2,y)=D_{x_2}K(x_1,x_2,y)$.

{\bf 4.} $K(x_1,0,y)=\frac{1}{y-x_1}$.

{\bf Proof.} These properties of $K(x_1,x_2,y)$ follow from definition, see the discussion in Section 5.1. The proof of uniqueness of the solution of the differential equation $(D_{x_1}-D_{x_2})K=0$ with the properties {\bf 1}, {\bf 2}, {\bf 4} is omitted. $\Box$

Based on this explicit series representation for the kernel $K(x_1,x_2,y)$ one can derive the following holonomic system of differential equations.

{\bf Theorem 5.3.2.} The function $K(x_1,x_2,y)$ satisfies the following differential equations:
$$D_{x_1}K-D_{x_2}K=0,$$
$$D_y^*K-D_{x_1}K=(r_1-1)(r_2-1)F\Big(r_1,r_2,s_1,\frac{x_1x_2}{t}\Big)$$
$$LK(x_1,x_2,y)=\frac{(s_1-1)t}{x_1x_2(x_1-y)(x_2-y)}F\Big(r_1-1,r_2-1,s_1-1,\frac{x_1x_2}{t}\Big)$$
where
$$D_y^*=\frac{d^2}{dy^2}\cdot y(y-1)(y-t)-\frac{d}{dy}\cdot y(y-1)(y-t)\Big(\frac{s_1}{y}+ \frac{s_2}{y-1}+\frac{s_3}{y-t}\Big)+r_1r_2y+\lambda$$
is the differential operator conjugate to $D_y$ and
$$L=\frac{(x_1-1)(x_1-t)}{(x_1-x_2)(x_1-y)}\cdot \frac{d}{dx_1}+\frac{(x_2-1)(x_2-t)}{(x_2-x_1)(x_2-y)}\cdot \frac{d}{dx_2}+\frac{(y-1)(y-t)}{(y-x_1)(y-x_2)}\cdot \frac{d}{dy}-$$
$$\frac{(r_1+r_2-2)x_1x_2y+(s_2+1-r_1-r_2)x_1x_2-(s_1+s_2-2)tx_1x_2+(s_1-1)t(x_1+x_2-y)}{x_1x_2(x_1-y)(x_2-y)}.$$

{\bf Theorem 5.3.3.} The kernel $K(x_1,x_2,y)$ admits the following integral representation
$$K(x_1,x_2,y)=\frac{s_1-1}{1-t}\Bigg(\frac{2(x_1-1)(x_2-1)}{1-t}\Bigg)^{1-s_2}\Bigg(\frac{2(x_1-t)(x_2-t)}{t(t-1)}\Bigg)^{1-s_3}\times$$
$$\int_0^1F(r_1-1,r_2-1,s_1-1,uq)q^{s_1-2}Q^{-1/2}\Big(-uq-v+w+1+Q^{1/2}\Big)^{s_2-1}\Big(-uq+v-w+1-Q^{1/2}\Big)^{s_3-1}dq$$

where
$$u=\frac{x_1x_2}{t},~v=\frac{(x_1-1)(x_2-1)(y-1)}{(t-1)^2},~w=\frac{(x_1-t)(x_2-t)(y-t)}{t(t-1)^2},$$
$$Q=(v-w)^2-2(v+w)(uq-1)+(uq-1)^2,$$
and we assume that $Q^{1/2} \sim v-w$ for $y \to \infty$.

After substitution the integral representation of Gauss hypergeometric function and the change of variables
\begin{equation} \label{q1q2}
q_1=\frac{-uq-v+w+1+Q^{1/2}}{2(1-pqu)},~~~q_2=\frac{-uq+v-w+1-Q^{1/2}}{2(1-pqu)}
\end{equation}
we obtain another integral formula for the kernel
$$K(x_1,x_2,y)=\Bigg(\frac{x_1x_2}{t}\Bigg)^{1-s_1}\Bigg(\frac{(x_1-1)(x_2-1)}{1-t}\Bigg)^{1-s_2}\Bigg(\frac{(x_1-t)(x_2-t)}{t(t-1)}\Bigg)^{1-s_3}\times$$
$$\frac{\Gamma(s_1)}{\Gamma(r_1-1)\Gamma(s_1-r_1)(1-t)}\int_{D}q_1^{s_2-1}q_2^{s_3-1}(1-q_1-q_2)^{s_1-r_1-1}\Bigg(1+\frac{v}{q_1}+\frac{w}{q_2}\Bigg)^{r_1-2}\frac{dq_1}{q_1}\cdot\frac{dq_2}{q_2}$$
where $D=\{(q_1,q_2),~0\leq p,q\leq 1\}$ and $q_1,q_2$ are parameterized by (\ref{q1q2}).

{\bf Theorem 5.3.4.} Let $s_1=r_1=1$. In this case the kernel $K(x_1,x_2,y)$ is given by
$$K(x_1,x_2,y)=\Bigg(\frac{2t-t(x_1+x_2+y)+x_1x_2y+tyP^{1/2}}{2t(1-x_1)(1-x_2)}\Bigg)^{s_2-1}\times$$
$$\Bigg(\frac{2t^2-t(x_1+x_2+y)+x_1x_2y+tyP^{1/2}}{2(t-x_1)(t-x_2)}\Bigg)^{s_3-1}\frac{1}{yP^{1/2}}$$
where
$$P=1-\frac{2x_1}{y}-\frac{2x_2}{y}+\frac{x_1^2}{y^2}+\frac{x_2^2}{y^2}-\frac{2x_1^2x_2}{ty}-\frac{2x_1x_2^2}{ty}+\frac{x_1^2x_2^2}{t^2}+\frac{2(2ty-y^2-t+2y)x_1x_2}{ty^2}.$$

{\bf Proof.} Set $r_1=1$ in the integral formula for the kernel above. After that write $(s_1-1)q^{s_1-2}dq=d(q^{s_1-1})$, integrate by parts and set $s_1=1$.  $\Box$

\subsection{The case of second order differential operators with more than 4 regular singular points}

Fix a natural number $n\geq 1$, pairwise distinct points $t_1,...,t_n\in\C$ such that $t_i\ne 0,1$ for $i=1,...,n$ and parameters $s_1,...,s_{n+2},r_1,r_2\in\C$ such that 
$$s_1+...+s_{n+2}=r_1+r_2+1.$$
Let
$$D_x=x(x-1)(x-t_1)...(x-t_n)\Bigg(\frac{d^2}{dx^2}+\Big(\frac{s_1}{x}+\frac{s_2}{x-1}+\frac{s_3}{x-t_1}+...+\frac{s_{n+2}}{x-t_n}\Big)\frac{d}{dx}\Bigg)+$$
$$+\lambda_1+\lambda_2x+...+\lambda_nx^{n-1}+r_1r_2x^n.$$
This is the most general second order differential operator with regular singularities at $x=0,1,t_1,...,t_n,\infty$ and having analytic solutions
 near $x=0,1,t_1,...,t_n$. 
 
 Let $f_{\vec{\lambda}}(x)$ be the unique solution of the equation $D_xf_{\vec{\lambda}}(x)=0$ analytic at $x=0$ and such that $f_{\vec{\lambda}}(x)=1+O(x)$. Write
$$f_{\vec{\lambda}}(x)=\sum_{i=0}^{\infty}P_ix_i$$ 
where $P_i,~i=0,1,...$ are polynomials in $\lambda_1,...,\lambda_n$ and $P_0=0$. One can show that the products $\{P_{i_1}...P_{i_n},~0\leq i_1\leq...\leq i_n\}$ form a basis in the vector space $\C[\lambda_1,...,\lambda_n]$.

Define structure constants $C_{i_1...i_{n+1}}^{j_1...j_n}$ by (\ref{strcm}).

Define the kernel $K(x_1,...,x_{n+1},y_1,...,y_n)$ as the generating function
by (\ref{genKm}).

{\bf Lemma 5.4.1.} The kernel $K(x_1,...,x_{n+1},y_1,...,y_n)$ is the unique function characterized by the following properties:

{\bf 1.} It is symmetric with respect to $x_1,...,x_{n+1}$ and with respect to $y_1,...,y_n$.

{\bf 2.} It has Laurent series expansion by non-negative powers of $x_1,...,x_{n+1}$ and by negative powers of $y_1,...,y_n$.

{\bf 3.} The following differential equation holds
$$\sum_{i=1}^{n+1} \frac{1}{(x_1-x_i)...\hat{i}...(x_{n+1}-x_i)}D_{x_i}K(x_1,...,x_{n+1},y_1,...,y_n)=0$$

{\bf 4.} $\displaystyle K(x_1,...,x_n,0,y_1,...,y_n)=\frac{1}{n!}\sum_{\sigma\in S_n}\frac{1}{(y_{\sigma_1}-x_1)...(y_{\sigma_n}-x_n)}.$

{\bf Theorem 5.4.1.} The following equation holds
$$K(x_1,...,x_{n+1},y_1,...,y_n)=\sum_{i_1,...,i_n\geq 0}F\Big(r_1+i_1+...+i_n,r_2+i_1+...+i_n,s_1,\frac{x_1...x_{n+1}}{t_1...t_n}\Big)\times$$
$$u_0^{i_1+...+i_n}P_{i_1,...,i_n}\Big(\frac{u_1}{u_0},...,\frac{u_n}{u_0}\Big)Q_{i_1,...,i_n}(y_1,...,y_n)$$
where $F$ is the Gauss hypergeometric function,
$$P_{i_1,...,i_n}(v_1,...,v_n)=$$
$$\sum_{j_1,...,j_n\geq 0}\frac{\prod\limits_{l=1}^{j_1+...+j_n}(s_2+i_1+...+i_n-l)\prod\limits_{l=1}^{j_1}(i_1-l+1)...\prod\limits_{l=1}^{j_n}(i_n-l+1)}{\prod\limits_{l=1}^{j_1}(s_3+l-1)...\prod\limits_{l=1}^{j_n}(s_{n+2}+l-1)}\frac{v_1^{j_1}...v_n^{j_n}}{j_1!...j_n!}$$
$$u_i=\frac{(x_1-t_i)...(x_{n+1}-t_i)}{t_i(t_i-1)(t_1-t_i)...\hat{i}...(t_n-t_i)},~i=1,...,n,~u_0=\frac{(x_1-1)...(x_{n+1}-1)}{(t_1-1)...(t_n-1)}.$$
Note that $u_0^{i_1+...+i_n}P_{i_1,...,i_n}\Big(\frac{u_1}{u_0},...,\frac{u_n}{u_0}\Big)$ are polynomials in $u_0,u_1,...,u_n$ for 
non-negative $i_1,...,i_n\in\Z$ and the sum in the definition of $P_{i_1,...,i_n}$ is finite in this case.

The functions $Q_{i_1,...,i_n}(y_1,...,y_n)$ are determined by the system of equations
$$\sum_{i_1,...,i_n\geq 0}\tilde{u}_0^{i_1+...+i_n}P_{i_1,...,i_n}\Big(\frac{\tilde{u}_1}{\tilde{u}_0},...,\frac{\tilde{u}_n}{\tilde{u}_0}\Big)Q_{i_1,...,i_n}(y_1,...,y_n)=\frac{1}{n!}\sum_{\sigma\in S_n}\frac{1}{(y_{\sigma_1}-x_1)...(y_{\sigma_n}-x_n)}$$
where $\tilde{u}_i=u_i|_{x_{n+1}=0}$.

We expect that the holonomic $D$-module generated by $K$ maps epimorphically at the generic point to the one described in Section 3.4 with a possible change of the cyclic 
vector. 

\subsection{The case of third order differential operators with 3 regular singular points}

Let 
$$D_x=x^2(x-1)^2\frac{d^3}{dx^3}+x(x-1)(a_1+a_2x)\frac{d^2}{dx^2}+(a_3+a_4x+a_5x^2)\frac{d}{dx}+a_6x+\lambda$$
where $a_1,...,a_6,\lambda$ are parameters.
This is the most general third order differential operator with regular singularities at $x=0,1,\infty$ and with analytic solutions near $x=0,1$.

There exists a unique solution $f_{\lambda}(x)$ of the differential equation 
$$D_xf_{\lambda}(x)=0$$
such that $f_{\lambda}(x)$ is analytic near $x=0$ and $f_{\lambda}(0)=1$. We have
$$f_{\lambda}(x)=\sum_{i=0}^{\infty}P_i(\lambda)x^i$$
where $P_i(\lambda)$ are polynomials in $\lambda$ of degree $i$ and $P_0(\lambda)=1$.

Since $P_i(\lambda)$, $i=0,1,...$ is a basis of the vector space $\C[\lambda]$, we can define structure constants of polynomial multiplication in this basis by (\ref{strc}).

Define kernel $K(x_1,x_2,y)$ as the generating function of these structure constants by (\ref{genK})

Introduce new parameters $b_1,b_2,c_1,c_2,c_3$ by
$$b_1+b_2=-a_1-1,~b_1b_2=a_3,~c_1+c_2+c_3=a_2-3,$$
$$c_1c_2+c_1c_3+c_2c_3=a_5-a_2+2,~c_1c_2c_3=a_6.$$

{\bf Lemma 5.5.1.} The kernel $K(x_1,x_2,y)$ satisfies the following differential equations
$$D_{x_1}K=D_{x_2}K,$$
$$D_y^*-D_{x_1}K=(c_1-1)(c_2-1)(c_3-1)\sum_{i=0}^{\infty}\prod_{l=0}^{i-1}\frac{(c_1+l)(c_2+l)(c_3+l)}{(b_1+l)(b_2+l)}\cdot \frac{(x_1x_2)^i}{i!}$$
where $D_y^*$ is the conjugate differential operator given by
$$D_y^*=-\frac{d^3}{dy^3}\cdot y^2(y-1)^2+\frac{d^2}{dy^2}\cdot y(y-1)(a_1+a_2y)-\frac{d}{dy}\cdot (a_3+a_4y+a_5y^2)+a_6y+\lambda$$

{\bf Theorem 5.5.1.} The kernel $K(x_1,x_2,y)$ is given by
$$K(x_1,x_2,y)=\frac{1}{y-1}\sum_{\substack{0\leq j\leq k\leq i+j,\\0\leq i}}\frac{(-1)^jk!}{j!(k-j)!(i+j-k)!}\cdot \frac{(x_1x_2)^i\Big((x_1-1)(x_2-1)\Big)^j}{(y-1)^k}\times$$
$$\prod\limits_{l=j+1}^k(b_1b_2+c_1c_2+c_1c_3+c_2c_3+a_4+l(c_1+c_2+c_3)+(1-l)(b_1+b_2)+l^2-l+1)\times$$
$$\frac{\prod\limits_{l=k}^{i+j-1}(c_1+l)(c_2+l)(c_3+l)}{\prod\limits_{l=0}^{i-1}(b_1+l)(b_2+l)}.$$

{\bf Remark 5.5.1.} Let $D_x$ be third order differential operator with symbol $(x-t_1)^2(x-t_2)^2(x-t_3)^2\frac{\partial^3}{\partial x^3}$, with regular singular points 
at $x=t_1,t_2,t_3$, and solutions near these points of the form $f(x)=(x-t_i)^{b_{i,j}}(1+O(x-t_i))$ for $i,j=1,2,3$. Here $t_i,b_{i,j}$ are generic parameters such that 
$\sum_{1\leq i,j\leq 3}b_{i,j}=3.$ Notice that any two differential operators with these properties are differ by an arbitrary constant. 

Define a function $K_3(x,y,z)$ as
$$K_3=F\Bigg(q_1\frac{(x-t_1)(y-t_1)(z-t_1)}{(x-t_3)(y-t_3)(z-t_3)},q_2\frac{(x-t_2)(y-t_2)(z-t_2)}{(x-t_3)(y-t_3)(z-t_3)}\Bigg)$$
where $q_1=\frac{(t_2-t_3)^3}{(t_1-t_2)^3},~q_2=\frac{(t_1-t_3)^3}{(t_2-t_1)^3}$ and $F$ satisfies the system of differential equations
$$\frac{\partial^3G}{\partial u_{11}\partial u_{12}\partial u_{13}}=\frac{\partial^3G}{\partial u_{21}\partial u_{22}\partial u_{23}}=\frac{\partial^3G}{\partial u_{31}\partial u_{32}\partial u_{33}}.$$
Here 
$$G=F\Bigg(\frac{u_{11}u_{12}u_{13}}{u_{31}u_{32}u_{33}},\frac{u_{21}u_{22}u_{23}}{u_{31}u_{32}u_{33}}\Bigg)\cdot\prod_{1\leq i,j\leq 3}u_{i,j}^{-b_{i,j}}.$$
Then $K_3(x,y,z)$ satisfies the differential equations
$$\frac{\partial K_3}{\partial x}=\frac{\partial K_3}{\partial y}=\frac{\partial K_3}{\partial z}.$$
It will be interesting to understand if $K_3$ is somehow connected with the kernel $K$ constructed in Theorem 5.5.1. It looks feasible that $K_3$ also satisfies to some kind of associativity condition. The kernel $K_3$ is similar to kernels constructed in Remarks 3.4.2 for $l=n$, and 3.3.2.

\addcontentsline{toc}{section}{Acknowledgements}

\section*{Acknowledgements}

We thank V. Golyshev, V. Rubtsov, D. van Straten for the discussion at early stages of this project. We are also grateful to V.Drinfeld, P.Etingof, D.Kazhdan and 
other participants of Geometric Langlands Seminar for their remarks and comments. 
 A.O. is grateful to IHES for invitations and excellent working atmosphere.

\end{document}